\documentclass[11pt]{article}

\pdfoutput=1

\usepackage{OC_energy_manuscript}

\makeatletter
\def\BState{\State\hskip-\ALG@thistlm}
\makeatother

\title{Lagrangian Relaxation for Continuous-Time Optimal Control of Coupled Hydrothermal Power Systems Including Storage Capacity and a Cascade of Hydropower Systems with Time Delays}

\author[1]{Chiheb Ben Hammouda\,\orcidlink{0000-0002-8386-0406}}
\author[2]{Eliza Rezvanova\,\thanks{eliza.rezvanova@kaust.edu.sa}\orcidlink{0000-0002-7039-2527}}
\author[2]{Erik von~Schwerin\,\orcidlink{0000-0002-2964-7225}}
\author[2,3,4]{Ra\'ul~Tempone\,\orcidlink{0000-0003-1967-4446}}
\affil[1]{Mathematical Institute, Utrecht University, Utrecht, The Netherlands}
\affil[2]{King Abdullah University of Science and Technology (KAUST), Computer, Electrical and Mathematical Sciences \& Engineering Division (CEMSE), Thuwal, Saudi Arabia}
\affil[3]{Chair of Mathematics for Uncertainty Quantification, RWTH Aachen University, Aachen, Germany.}
\affil[4]{Alexander von Humboldt Professor in Mathematics for Uncertainty Quantification, RWTH Aachen University, Aachen  Germany.}

\begin{document}
	\date{}
\maketitle

\begin{abstract}
	This work considers a short-term, continuous time setting characterized by a coupled power supply system controlled exclusively by a single provider and comprising a cascade of hydropower systems (dams), fossil fuel power stations, and a storage capacity modeled by a single large battery.

Cascaded hydropower generators introduce time-delay effects in the state dynamics, which are modeled with differential equations, making it impossible to use classical dynamic programming. We address this issue by introducing a novel Lagrangian relaxation technique over continuous-time constraints, constructing a nearly optimal policy efficiently. This approach yields a convex, nonsmooth optimization dual problem to recover the optimal Lagrangian multipliers, which is numerically solved using a limited memory bundle method. At each step of the dual optimization, we need to solve an optimization subproblem. Given the current values of the Lagrangian multipliers, the time delays are no longer active, and we can solve a corresponding nonlinear Hamilton--Jacobi--Bellman (HJB) Partial Differential Equation (PDE) for the optimization subproblem. The HJB PDE solver provides both the current value of the dual function and its subgradient, and is trivially parallelizable over the state space for each time step.

To handle the infinite-dimensional nature of the Lagrange multipliers, we design an adaptive refinement strategy to control the duality gap. Furthermore, we use a penalization technique for the constructed admissible primal solution to smooth the controls while achieving a sufficiently small duality gap. Numerical results based on the Uruguayan power system demonstrate the efficiency of the proposed mathematical models and numerical approach.

\textbf{Keywords} Renewable energy, power supply systems, short-term continuous control, Lagrangian relaxation, cascade of hydropower generators with time delays, optimal control, dynamic programming algorithm, limited memory bundle method, duality gap

\end{abstract}

\thispagestyle{plain}

\section{Introduction}
Considering the global large-scale  development of renewable energy resources and the high uncertainty and variability of these sources, the control optimization of energy supply systems, particularly hydropower systems, has become increasingly crucial for the efficient use of production facilities to balance the instantaneous power demand. In this work, we propose a novel mathematical modeling and numerical framework for the optimal management of large-scale power systems in a short-term planning horizon, where the system comprises a cascade of hydropower stations with time delays, fossil fuel power stations (FFSs), and a storage capacity modeled by a single battery. To do so, we pose a time-continuous optimal control (OC) problem and design a primal-dual strategy by introducing a Lagrangian relaxation technique over continuous-time constraints.

Some previous studies~\cite{glanzmann2005supervisory,setz2008application,korobeinikov2010optimizing,ribeiro2012optimal,paredes2014huc, hamann2017mpc,mathur2020huc} have considered only the hydropower system, whereas others have studied the coupled hydrothermal~\cite{xi1999scheduling} or hydro-wind~\cite{hug2011predictive} power systems. These studies applied discrete-time discrete-space models to create hourly scheduling, leading to nondeterministic polynomial-time hard mixed-integer programming problems, which are usually challenging to solve. In particular, several studies~\cite{glanzmann2005supervisory, setz2008application,hug2011predictive, hamann2017mpc} have investigated a cascade of run-of-the-river dams to determine an optimal schedule using model predictive control (MPC).  
Compared to the present work, they modeled river dynamics instead of dam dynamics, using the water level rather than the dam volume as the state-space variable. The resulting OC problems within the model predictive control framework are solved using discrete optimization methods, including linear and nonlinear programming methods, depending on the objective function \cite{korobeinikov2010optimizing,glanzmann2005supervisory,setz2008application,hug2011predictive,hamann2017mpc}, (mixed) integer programming \cite{xi1999scheduling,paredes2014huc,mathur2020huc}, and nonlinear network flow algorithms \cite{xi1999scheduling}.

The authors of \cite{ferreira2015continuous} presented a continuous-time formulation of the OC of a cascade of hydropower stations without accounting for time delays. Moreover, to solve the corresponding OC problem, they used a different approach based on a formulation of an infinite-dimensional optimization problem with a quadratic functional subject to cone constraints, and they derived sufficient conditions for local minimizers. This approach is problem-specific and difficult to generalize to a broader class of (coupled) power systems, particularly those accounting for hydrosystem time delays. One important aspect that we address in this work is the time delay for the water flow in the dynamics of a cascade of dams. Some research \cite{nilsson1996mixed,guan1997optimization, xi1999scheduling,hamann2017mpc} has incorporated delay constraints but used a different approach than that presented here, based on a discrete-time framework and mixed-integer programming optimization. These studies decompose the hydrothermal problem into thermal and hydro-subproblems and further decompose the hydro-subproblem into several subproblems, usually mixing a continuous-variable optimization problem determining the generation levels of all units in the entire river and pure integer problems determining the hydro commitment states, one for each unit. In \cite{Zhang@formationcontrol}, the authors employ a different approach to solving the OC problem: instead of associating the HJB equation to the value function, reinforcement learning with fuzzy logic is used to minimize the cost function. Some other works study multiagent systems \cite{Deng2023distribalgogames, Deng2024distribalgo} and use distributed algorithms for nonsmooth optimization, which is also different from our approach. The cost functions of the agents in these works are revealed to the players only after the decisions are made. Furthermore, their problem settings do not contain time-delays in the system dynamics. In this work, we study a centralized decision maker, which gives us an opportunity to exploit our knowledge of the system structure and the cost functions characterizing all power sources to our advantage.

In this work, we design a primal-dual strategy to formulate and solve the corresponding continuous-time OC problem, leading to efficient policies minimizing the cost associated with power production subject to demand constraints and the dynamics and constraints of the hydropower and battery system. Because of the time delays in the system dynamics, the primal problem cannot be solved directly by the continuous dynamic programming principle. To overcome this undesirable feature, we pose the corresponding dual problem by introducing time-continuous Lagrangian relaxation to penalize the time-delay constraints in the system dynamics. The resulting dual-cost function, together with the appropriate dynamics and constraints, creates a standard continuous OC problem. 
Based on the dual-problem formulation, we iteratively optimize the approximated dual function with respect to (w.r.t.) the Lagrangian multipliers. Each iteration involves the subproblem of solving the associated Hamilton--Jacobi--Bellman (HJB) partial differential equation (PDE) with an upwind finite-difference scheme~\cite{sun2015HJBscheme} adapted to this context.
From the subproblem, we obtain the dual objective function and its subgradient.

Because the resulting dual problem is nonsmooth and convex, we solve it using the limited memory bundle method (LMBM)~\cite{LMBM_PhD,haarala2007globally,karmitsa2012comparing,bagirov2014introduction}. Moreover, we design an adaptive refinement strategy for the Lagrangian multipliers to control the dual optimization error better. Finally, for the constructed admissible primal solution, we employ a penalization technique to smooth the controls while achieving a sufficiently small desired duality gap. The numerical experiments and results based on the Uruguayan power system demonstrate the efficiency of the proposed mathematical models and numerical approach. 
Preliminary work was conducted in the master theses in~\cite{caballero2019thesis} and~\cite{rezvanova2021thesis}.

Our contributions can be summarized as follows:
\begin{itemize}
	\item We propose a novel time-continuous OC framework for optimal management of large-scale power systems with time-delays. This system includes a cascade of hydropower stations, FFSs, and a storage capacity modeled by a single battery.	Compared to the discrete-time discrete-space formulation presented in the mentioned studies, the proposed continuous-time and continuous-space modeling offers an advantage of decoupling the process of model development from numerical approximation, thereby enhancing model fidelity. It is grounded in an underlying physical model that encapsulates the dynamic behavior of dam and battery systems, as well as hydropower generation processes. Thus, our formulation provides a continuous curve of control over time, allowing its application for any time stepping and eliminating the need for an ad hoc interpolations often employed in a discrete setting.	
	\item The time delays in the system are dealt with by utilizing a primal-dual approach, incorporating Lagrangian relaxation over continuous-time constraints to handle the complexities of the system dynamics.		
	\item We designed an iterative process for optimizing the dual function concerning Lagrangian multipliers. This process includes solving the HJB partial differential equation and employing the LMBM for solving the resulting nonsmooth and convex dual problem.	
	\item To achieve a sufficiently small desired duality gap, we designed an adaptive refinement of Lagrangian multipliers approximation strategy.	
	\item We designed a penalization technique to construct a smooth admissible primal solution, which can be implemented in practice.
	\item To demonstrate the efficiency of the proposed mathematical models and numerical approach we study the Uruguayan power system. The management and regulation of Uruguayan electricity market are overseen by the Administración del Mercado Eléctrico (ADME). Currently, the methodologies employed for power system scheduling in this sector utilize Stochastic Dual Dynamic Programming (SDDP), SimSEE, and Plexos \cite{adme2023report}. These methodologies are characterized as discrete-time, discrete-space approaches, using scenario-based risk assessment to address variability sources, and do not include the storage capacity in these scenarios. Moreover, the time delays inherent in the dams dynamics are not taken into account, which may lead to infeasible solutions.		
\end{itemize}

The outline of this paper is as follows. Section 2 introduces the model of the complete linked system, including a detailed description of the cost, power generation, and state dynamics for each power source. Next, Section 3 formulates the primal OC problem for the complete linked system with time delays. Subsequently, we propose a continuous-time Lagrangian relaxation technique to address the time delays in hydropower dynamics, resulting in a dual problem formulation. Section 3 also covers the process of obtaining nearly optimal admissible controls that are sufficiently smooth. Then, Section 4 delves into the numerical approach for solving the primal problem, following the method described in Section 3. We present the algorithm in detail, explaining its components. In particular, we discuss the numerical solution of the HJB equation, aspects of dual-problem optimization, the parameter choice for constructing the sufficiently smoothed primal solution from the dual solution, and various numerical errors arising from the proposed approach. Finally, Section 5 focuses on the Uruguayan power grid as an illustrative example, demonstrating the numerical results when solving the OC problem for managing the related coupled power system.
\section{Power System Model}\label{sec:Mathematical Modeling of Energy Systems}
We consider a setting characterized by a coupled power supply system controlled exclusively by a single provider (monopoly). The central operator controls both the supply and price of power to balance the instantaneous power demand. The efficient use of power-generation facilities is crucial for two reasons: i) generated profit ensures resources for infrastructure improvement, and ii) failure to meet the power demand introduces a very high cost. We introduce the demand and power balance equation in Section~\ref{sec:Demand modeling}.

Both the demand and available outputs of various power sources have some degree of uncertainty, even day-to-day, due to uncontrollable factors, such as the weather. This uncertainty can be included in the modeling by introducing stochastic dynamics. The present work is 
restricted to the deterministic case; however, the suggested approach was developed with stochastic dynamics in mind.

Mathematically, the control problem is challenging because the central operator must meet the demand at each time and for each realization of the stochastic dynamics, making the problem infinite-dimensional.  

We specifically consider a setting where the operator manages a hybrid electricity production system consisting of various power sources. The system of study comprises only controllable sources, where the production can be controlled over time, subject to constraints. The power sources include the following: 
\begin{itemize}
	\item a hydropower network consisting of $N_D$ dams, where some dams are serially connected along the same river so that the downstream dams depend on the outflow of upstream dams,
	\item $N_F$ FFSs, and
	\item a storage capacity modeled by a single large battery.
\end{itemize}
Section~\ref{sec:Cost and Power Generation Equations} introduces each power source's
cost, power generation, and state dynamics. The power generated by these sources is meant to satisfy the effective demand 
after subtracting the power generated by noncontrollable sources (e.g., wind and solar) from the target demand. In future work, we aim to add noncontrollable sources to the studied power systems, adding further stochasticity to the optimal problem and further challenges in modeling and numerics.

\subsection{Demand Model}\label{sec:Demand modeling}
We optimize the dispatch of controllable power sources, explained in Section~\ref{sec:Cost and Power Generation Equations}, to meet the effective demand, $D_E(t)$, at any time $t$. The effective demand is 
obtained by subtracting the power generated by noncontrollable sources from the total local demand, $D(t)$, and exports, $E(t)$. 
With a finite time horizon, $T$, the power generated by the power system must satisfy the effective demand given by 
\begin{align}\label{eq:demand_eq}
	D_E(t) & = D(t)+E(t)- P_{WS}(t), && 0 \le t\le T,
\end{align}
where $P_{WS}(t)$ is the generated wind and solar power.

We restrict the work to deterministic models for effective demand. 
In short-term (one-day) optimizations, reasonably accurate predictions for the future electricity demand are typically available to the operators. 

Given the effective demand and controllable resources, the instantaneous power balance equation is 
\begin{align}\label{eq:power_balance_equation}
	D_E(t)  & =\sum_{i=1}^{N_D} P_H^{(i)}(t)  +\sum_{i=1}^{N_F}P_F^{(i)}(t) +P_A(t), && 0 \le t \le  T,
\end{align}
where $P_H^{(i)}(t)$, $P_F^{(i)}(t)$, and $P_A(t)$ denote the instantaneous power supplied by the dams, FFSs, and battery, respectively (see Sections~\ref{sec:Dams Modeling}, \ref{sec:Fossil Fuel Stations (FFS) Modeling}, and~\ref{sec:Battery Modeling}).

\subsection{Cost, Power Generation, and State Dynamics}\label{sec:Cost and Power Generation Equations}
For mathematical convenience and to avoid needlessly poor conditioning due to varying scales of  different quantities in the numerical computations, the following model dynamics are normalized to take values in the range $[0,1]$. 
Where absolute maxima and minima are available, these are used in the normalization; 
otherwise, upper and lower bounds respecting the order of magnitude are used. 

\subsubsection{Hydropower Network Model}\label{sec:Dams Modeling}
\begin{SCfigure}
	\centering
	\caption{Example configuration of dams with the flow directed downward. \\
		$\bullet$ Each dam in the network is identified by a unique index, $i$.\\
		$\bullet$ Any dam $i$ that has another dam downstream is associated with
		a constant, $\tau_i$, representing the time it takes for the water to 
		reach the downstream neighbor.\\
		$\bullet$ Any dam $i$ is associated with a set, $\mathcal{B}(i)$, containing 
		the indices of its nearest upstream dams along any fork of the river.\\
		$\bullet$ No dam has more than one downstream neighbor.
	}
	\includegraphics[width=0.65\textwidth]{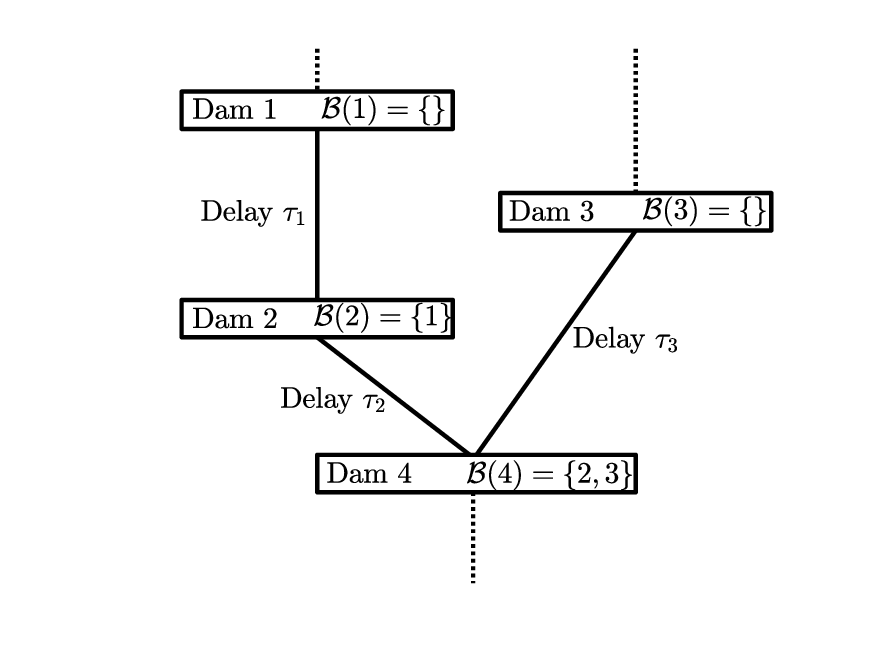}
	\label{fig:dam_notations}
\end{SCfigure}

Each dam is represented by the volume of water in its reservoir, $v(t)$. The instantaneous change in this volume is the difference between the inflow of water from rivers and rain, denoted by $I(t)$, and the water outflow. We assume two types of controls are related to the water outflow: i) the control of the flow through the turbines, denoted by $\phi_{Tur}$, to produce energy, and ii) the control of how much water is spilled, denoted by $\phi_S$. In the dynamics of the connected dams, the inflow to any downstream dam depends on the controlled outflow of the nearest upstream dams. The entire outflow of an upstream dam is assumed to reach the following dam with a constant time delay, $\tau>0$. The OC must respect this time delay in continuous time. Rivers can merge but not branch out again; thus, each dam has at most one downstream neighbor but can have several upstream neighbors. The set of upstream neighbors of a dam is denoted as $\mathcal{B}$ (see Figure~\ref{fig:dam_notations}).
By practical limitations, the volume of each reservoir is bounded from above and below. We 
define, for the $i$th dam the minimum and maximum values $\underline{v}^{(i)}$ and $\overline{v}^{(i)}$, respectively. The normalized volume $\hat{v}^{(i)}(t)=v^{(i)}(t)/\left(\bar{v}^{(i)}-\underline{v}^{(i)}\right)$
takes values in the range $[0,1]$.

\paragraph{Dam dynamics}

For the independent dams, the (dimensionless) dynamics\footnote{We write the ordinary differential equations in differential form, adhering to the convention for It\^{o} stochastic differential equations; see Remark~\ref{rem:v_SDE}.} of the $i$th dam are modeled by
\begin{align}\label{eq: dams_independent_dynamcis}
	d\hat{v}^{(i)}(t) &= \frac{1}{\bar{v}^{(i)} -  \underline{v}^{(i)}} \left(I^{(i)}(t)- \phi_{Tur}^{(i)}(t)- \phi_S^{(i)}(t)\right) dt, && 0 \le t \le  T,
\end{align}
with the initial condition $\hat{v}^{(i)}(0) =\hat{v}^{(i)}_0\in[0,1]$. 
The requirement that $0\le \hat{v}^{(i)}(t) \le 1$, for $t>0$, is enforced by imposing the following outflow constraints on the control:
\begin{align}\label{eq: dams_independent_constraints}
	& \begin{cases}
		\phi_{Tur}^{(i)}(t) + \phi_S^{(i)}(t) \geq I^{(i)}(t),  & \text{if $\hat{v}^{(i)} (t)= 1$,} \\
		\phi_{Tur}^{(i)}(t) + \phi_S^{(i)}(t) \leq I^{(i)}(t),  & \text{if $\hat{v}^{(i)} (t)= 0$.} 
	\end{cases}
\end{align}

Naturally, the outflow constraints cannot be satisfied for arbitrary inflows, $I^{(i)}$, 
because the maximal controlled outflow from a dam is bounded. 
This study assumes that the inflows are such that the constraints can be satisfied.

For a downstream dam, the dynamics and constraints are adjusted by the outflow of the 
neighboring upstream dams. Then,~\eqref{eq: dams_independent_dynamcis} is replaced by
\begin{align}\label{eq: dams_connected_dynamcis}
	d\hat{v}^{(i)}(t) = \frac{1}{\bar{v}^{(i)} -  \underline{v}^{(i)}} \left(I^{(i)}(t)+\sum_{j\in\mathcal{B}(i)}\left(\phi_{Tur}^{(j)}(t-\tau_j)+ \phi_S^{(j)}(t-\tau_j)\right)- \phi_{Tur}^{(i)}(t)- \phi_S^{(i)}(t)\right) dt, 
\end{align}
for $0 \le t \le T$,
and~\eqref{eq: dams_independent_constraints} is replaced by
\begin{align}\label{eq: dams_connected_constraints}
	& \begin{cases}
		\phi_{Tur}^{(i)}(t) + \phi_S^{(i)}(t) \geq 
		I^{(i)}(t)+\sum_{j\in\mathcal{B}(i)}\left(\phi_{Tur}^{(j)}(t-\tau_j)+ \phi_S^{(j)}(t-\tau_j)\right),  
		& \text{if $\hat{v}^{(i)} (t)= 1$,}\\
		\phi_{Tur}^{(i)}(t) + \phi_S^{(i)}(t) \leq
		I^{(i)}(t)+\sum_{j\in\mathcal{B}(i)}\left(\phi_{Tur}^{(j)}(t-\tau_j)+ \phi_S^{(j)}(t-\tau_j)\right),  
		& \text{if $\hat{v}^{(i)}(t) = 0$.}\\
	\end{cases}
\end{align}
The dynamics and constraints in~\eqref{eq: dams_independent_dynamcis} and~\eqref{eq: dams_independent_constraints} are special cases of~\eqref{eq: dams_connected_dynamcis} and~\eqref{eq: dams_connected_constraints}. We include the known upstream turbine flows and spillages at times $t<0$ in the given inflow model $I^{(i)}$ and  
define the controls to satisfy $\phi_{Tur}^{(j)}(t)\equiv\phi_{S}^{(j)}(t)\equiv 0$, for $t<0$.

The flow through the turbines and water spillage are both bounded, $0\leq \phi_{Tur}^{(i)}\leq \bar{\phi}_{Tur}^{(i)}(\hat{v}^{(i)})$ and $0\leq \phi_S^{(i)}\leq \bar{\phi}_S^{(i)}(\hat{v}^{(i)})$, where the upper bounds 
on the flow may depend on the volume (as a function of the height) of the water in the reservoir. Given these constraints, we introduce the normalized controls: 
\begin{align}
	\label{eq:scaling_dams_controls}
	\hat{\phi}_{Tur}^{(i)}(t) &= \frac{\phi_{Tur}^{(i)} (t)}{\bar{\phi}_{Tur}^{(i)}\left(\hat{v}^{(i)}(t)\right)}, 
	& 
	\hat{\phi}_S^{(i)}(t) = \frac{\phi_S^{(i)} (t)}{\bar{\phi}_S^{(i)}\left(\hat{v}^{(i)}(t)\right)},
	&& \text{for $1 \le i \le N_D$,}
\end{align}
where $\bar{\phi}_{Tur}^{(i)}(\hat{v}^{(i)})$ and $\bar{\phi}_S^{(i)}(\hat{v}^{(i)})$ depend on the dam architecture and hydroturbine technology. 
In~\eqref{eq:model_phimax}, we illustrate examples of $\bar{\phi}^{(i)}_T(\hat{v}^{(i)})$ of the dams considered in these experiments. Moreover, for ease of presentation,  $\bar{\phi}_{Tur}^{(i)}(\hat{v}^{(i)})$ and $\bar{\phi}_S^{(i)}(\hat{v}^{(i)})$ are denoted as $\bar{\phi}_{Tur}^{(i)}$ and $\bar{\phi}_S^{(i)}$, respectively.

\paragraph{Hydropower generation}

The power generated by a dam primarily depends on two quantities: i) the flow $\phi_{Tur}$ through the turbines, and ii) the difference between the upstream water level, $H$, and the downstream water level, $h_0$, measured relative to a common reference frame (see, e.g., Chapter 6.2.2 of \cite{gharehpetian2017generationsystems}, \cite{ paish2002hydropower, purohit2008hydropowerindia}). Increasing the flow $\phi_{Tur}+\phi_S$ has the side effect of slightly decreasing the height difference. The strength of this effect for different flows and water levels depends on the dam geometry. Motivated by the Uruguayan power system data used in Section~\ref{sec:num_experiments}, we consider this effect by linearly modeling the generated power for the $i$th dam, $P_H^{(i)}$, in $\hat{\phi}^{(i)}_{Tur}$ as
\begin{subequations}\label{eq:HP_linearized model}
	\begin{align}
		P_H^{(i)}(\hat{\phi}^{(i)}_{Tur}(t),\hat{v}^{(i)}) & = S^{(i)}(\hat{v}^{(i)}) \hat{\phi}^{(i)}_{Tur}(t), \\
		S^{(i)}(\hat{v}^{(i)}) &= \eta^{(i)} \overline{\phi}^{(i)}_{Tur}(\hat{v}^{(i)}) \Big( H^{(i)}(\hat{v}^{(i)}) - h_{0}^{(i)} - d^{(i)} \overline{\phi}^{(i)}_{Tur} \Big),
	\end{align}
\end{subequations}
where $h_0^{(i)}$ is considered independent of $\phi_{Tur}+\phi_S$, and the dam specific, positive model parameters $\eta^{(i)}$ and $d^{(i)}$ are calibrated to the data.

The upstream water level, $H^{(i)}$, may have a complicated dependency on the volume, $\hat{v}^{(i)}$. The water level is typically modeled as a polynomial of $\hat{v}^{(i)}$ (see~\cite{ribeiro2012optimal} and \eqref{eq:hight_volume_dependence}). These polynomials should be increasing and satisfy the constraints implied by the corresponding constraints on the volumes. 

Finally, we model the total cost of generating power by the $i$th dam in $[0, T]$ using 
\begin{align}\label{eq:dams_power_cost_eq}
	C_H^{(i)}(T)&= K_H^{(i)} \int_{0}^{T} \left(\phi_{Tur}^{(i)} (t)+\phi_S^{(i)}(t) \right) dt,
\end{align}
where $K_H^{(i)}>0$ is the given water cost of the $i$th dam.
\begin{remark}[On quadratic hydropower generation models]
	A frequently used alternative to~\eqref{eq:HP_linearized model} is
	\begin{equation}\label{eq:dams_power_eq}
		P_H^{(i)}\left( \phi_{Tur}^{(i)},\phi_S^{(i)};v^{(i)}\right)= \eta^{(i)} \; \phi_{Tur}^{(i)} \left(H^{(i)} (v^{(i)})-h_0^{(i)}-d^{(i)} \left( \phi_{Tur}^{(i)}+\phi_S^{(i)}\right)\right).
	\end{equation}
	This is quadratic in the controls $\phi^{(i)}_{Tur}$ and $\phi^{(i)}_S$; thus,
	the power balance constraint changes from linear to nonlinear, requiring 
	adjustments to the numerical minimization of the Hamiltonian.
\end{remark}
\begin{remark}[The cost of water]
	The coefficients $\{K_H^{(i)}\}_{i=1}^{N_D}$ in \eqref{eq:dams_power_cost_eq} reflect the future value of water, that is, how valuable the hydropower sources (dams) are compared to other sources. These coefficients represent future profits or losses depending on the ability 
	of the dams to substitute other costly sources, such as fuel. In practice, $\{K_H^{(i)}\}_{i=1}^{N_D}$ are outputs of mid-term optimization models, and in this work, we consider them  inputs.	
\end{remark}
\begin{remark}[Reversible turbines]
	In modern reversible hydroelectric power stations, the turbine can be reversed to pump water from a downstream to an upstream reservoir. Water is usually pumped upstream at times of low electricity demand to build up reserves to produce energy during peak hours, thus balancing the load and making a profit on the price difference. We leave adding this feature to the proposed dam modeling and studying its consequences on the SOC problem for future work.
\end{remark}

\subsubsection{Storage Capacity Model}\label{sec:Battery Modeling}
In this power network system, we add storage capacity modeled as a single large battery, which is fully controllable, and its state is its capacity, $a(t)$. There are many battery models of varying complexity and perspective developed for different application areas \cite{tomasov2019battery, tamilselvi2021batteryreview}. In this work, we use a simple mathematical model presented below in (\ref{eq: battery_dynamcis}), based on a circuit model of the state of charge of a lithium-ion battery (see eq. (9a) in \cite{brucker2021batterymodel}). We assume that it has no energy loss when used and that no associated operational costs or aging effects exist. 

We let $\bar{A}$ be the maximum battery capacity and $P_A(t)$ be the instantaneous power supplied by the battery to the system at time $t$. 
As the battery can be charged, $P_A(t)<0$, and discharged, $P_A(t)>0$, we define 
the extremes, $\underline{P}_A$ and $\overline{P}_A$, s.t. $P_A(t) \in [\underline{P}_A, \:\overline{P}_A]$, $\forall t \in [0, T]$. Generally, $|\underline{P}_A| < \overline{P}_A$ (i.e., the battery can supply higher power than it can absorb from the 
grid).

As above, we use the normalized capacity, $\hat{a}(t)$, and normalized control, $\hat{\phi}_A(t)$, satisfying
$P_A(t) = \overline{P}_A \hat{\phi}_A(t)$. 
With this normalization, $\hat{\phi}_A(t) \in [- C_{\text{bat}} ,1]$, where the constant 
$C_{\text{bat}}=-\underline{P}_A/\overline{P}_A>0$ is a characteristic of the battery.
The normalized capacity dynamics, $\hat{a}(t)$, is
\begin{align}\label{eq: battery_dynamcis}
	d\hat{a}(t)&= 
	-\frac{\overline{P}_A}{\bar{A}} \hat{\phi}_A(t) dt, && 0 \le t \le  T,
\end{align}
with the initial condition $\hat{a}(0)=\hat{a}_0\in [0,1]$. The requirement that $0\le \hat{a}(t) \le 1$, for $t>0$, is enforced through the constraints on charging and discharging at the extremes:
\begin{align}\label{eq: battery_constraints}
	& \begin{cases}
		P_A(t) \ge 0, & \text{if $\hat{a}(t)= 1$,} \\
		P_A(t) \le 0, & \text{if $\hat{a}(t)= 0$.} 
	\end{cases}
\end{align}

\subsubsection{Fossil Fuel Station Model}
\label{sec:Fossil Fuel Stations (FFS) Modeling}

This work considers FFSs to be without a state but with a control, which is their power. Moreover, we assume that, for the $i$th FFS, both the cost associated with the generated power, $K^{(i)}_F>0$, and the maximum available power, $\bar{P}_F^{(i)}>0$, can change weekly or daily but are considered constant during at least one entire day.

We model the power generation of the $i$th station, $P_F^{(i)}(t)$, as linear w.r.t. the normalized fossil fuel control, $\hat{\phi}_F^{(i)}(t)$, that is,
\begin{align}\label{eq:FFS_power_eq}
	P_F^{(i)}\left(\hat{\phi}_F^{(i)}(t)\right)&= \bar{P}_F^{(i)} \hat{\phi}_F^{(i)}(t), && 0 \le t \le  T.
\end{align}
We ignore the finite start-up time and minimum stable capacity of FFSs and model the total cost of generating power using the $i$th FFS between time $t=0$ and time $T$, as follows:
\begin{align}\label{eq:FFS_power_cost_eq}
	C_F^{(i)}(T)  &= 
	K_F^{(i)}\bar{P}_F^{(i)} \int\limits^{T}_0 \hat{\phi}_F^{(i)}(t) dt.
\end{align}
\begin{remark}
	In the current work, we do not deal with the unit commitment optimization (nonlinear mixed-integer optimization) problem \cite{kazarlis1996genetic,bertsimas2012adaptive,saravanan2013solution}, which determines the optimal operation schedule of the different power units to meet the demand under various constraints (e.g., start-up and shut-down costs). We leave this consideration in the proposed models for future work.
\end{remark}

\subsection{Model Summary}
\label{sec:model_summary}

The state and control vectors are denoted by $\mathbf{X}(t)$ and $\hat{\boldsymbol{\phi}}(t)$, respectively: 
\begin{align}
	\mathbf{X}(t)              & := \left( \hat{v}^{(1)}(t),\dots,\hat{v}^{(N_D)}(t), \hat{a}(t) \right),\\
	\hat{\boldsymbol{\phi}}(t) & := \left( \hat{\phi}_{Tur}^{(1)}(t),\dots,\hat{\phi}_{Tur}^{(N_D)}(t), \hat{\phi}_S^{(1)}(t),\dots,\hat{\phi}_S^{(N_D)}(t), \hat{\phi}_F^{(1)}(t),\dots,\hat{\phi}_F^{(N_F)}(t), \hat{\phi}_A (t) \right).
\end{align}

For fixed initial conditions and deterministic forecasts of $D_E$ and $I^{(i)}$, 
the control is a measurable\footnote{We add extra smoothness requirements in 
	Section~\ref{sec:Admissible Controls Smoothing: Computation of Penalized Controls}} 
function $\hat{\boldsymbol{\phi}}:[0,T]\to[0,1]^{N_D+N_F}\times[-C_\text{bat},1]$, 
and the controlled state evolution is $\mathbf{X}(t):[0,T]\to[0,1]^{N_D+1}$.
Together, these satisfy the following dynamics and constraints. 
\begin{dynamics}[Resource Dynamics]\label{dynamics}
	The components of the state vector $\mathbf{X}(t)$ (i.e., the normalized volumes, 
	$\hat{v}^{(i)}$, for $i=1,\dots,N_D$, and battery capacity, $\hat{a}$) 
	satisfy the controlled dynamics (\ref{eq: dams_connected_dynamcis}) and (\ref{eq: battery_dynamcis}) in $[0,T]$,
	where $\phi_{Tur/S}^{(j)}(t) = \bar{\phi}_{Tur/S}^{(j)}\left(\hat{v}^{(j)}(t)\right)\hat{\phi}_{Tur/S}^{(j)}(t)$. All initial conditions, $\hat{v}^{(i)}(0)$ and $\hat{a}(0)$, are in $[0,1]$.
\end{dynamics}

\begin{constraints}[Instantaneous Control Constraints]\label{instantaneous_constraints}
	The ranges of the normalized controls are
	\begin{subequations}\label{eq:norm_con}
		\begin{align}
			0 & \leq \hat{\phi}_{Tur}^{(i)}(t), \hat{\phi}_S^{(i)}(t), \hat{\phi}_F^{(j)}(t) \leq 1, 
			&& i=1,2,\dots,N_D, j=1,2,\dots,N_F\\
			- C_{\text{bat}} & \leq \hat{\phi}_A(t) \leq 1, && ~
		\end{align}
	\end{subequations}
	with the further state-dependent constraints 
	\begin{subequations}\label{eq:con_a}
		\begin{align}
			\hat{\phi}_A(t) & \geq 0, && \text{if $\hat{a}(t)= 1$,}\\
			\hat{\phi}_A(t) & \leq 0, && \text{if $\hat{a}(t) = 0$,}
		\end{align}
	\end{subequations}
	and, for any $i=1,2,\dots,N_D$ s.t. $\mathcal{B}(i)$ is empty, (\ref{eq: dams_independent_constraints}) holds. 
	Furthermore, given the deterministic forecast of the effective demand, $D_E(t)$, the controls are constrained to satisfy
	\begin{align}
		D_E(t) & = 
		\sum_{i=1}^{N_D} S^{(i)}(\hat{v}^{(i)}) \hat{\phi}^{(i)}_{Tur}(t)+\sum_{i=1}^{N_F} \bar{P}_F^{(i)} \hat{\phi}_F^{(i)}(t)+ \overline{P}_A \hat{\phi}_A(t).
	\end{align}
\end{constraints}

\begin{constraints}[Time-Delayed Control Constraints]\label{time_delayed_constraints}
	For any $i=1,2,\dots,N_D$ s.t. $\mathcal{B}(i)$ is nonempty, (\ref{eq: dams_connected_constraints}) holds. 
\end{constraints}
As a consequence of the dynamics and constraints, $\mathbf{X}(t)\in [0,1]^{N_{D}+1}$, for all times.
\section{Optimal Control Formulation}\label{sec:SOC_formulation}
Section~\ref{sec:Statement of the Optimal Control problem of the Complete Linked System} formulates the OC problem of the complete linked system with the cost, power generation, and state dynamics described in Section~\ref{sec:Cost and Power Generation Equations}, given the 
demand~\eqref{eq:demand_eq}. Then, Section~\ref{sec:Continuous-Time Lagrangian Relaxation} introduces the continuous-time Lagrangian relaxation, leading to a dual-problem formulation to address the time delays in the hydropower dynamics. Finally, Sections~\ref{sec:Construction of the admissible Controls for the Primal Problem} and~\ref{sec:Admissible Controls Smoothing: Computation of Penalized Controls} explain how to obtain nearly optimal admissible controls that are sufficiently smooth.

\begin{figure}[h!]
	\centering
	\includegraphics[width=0.9\textwidth]{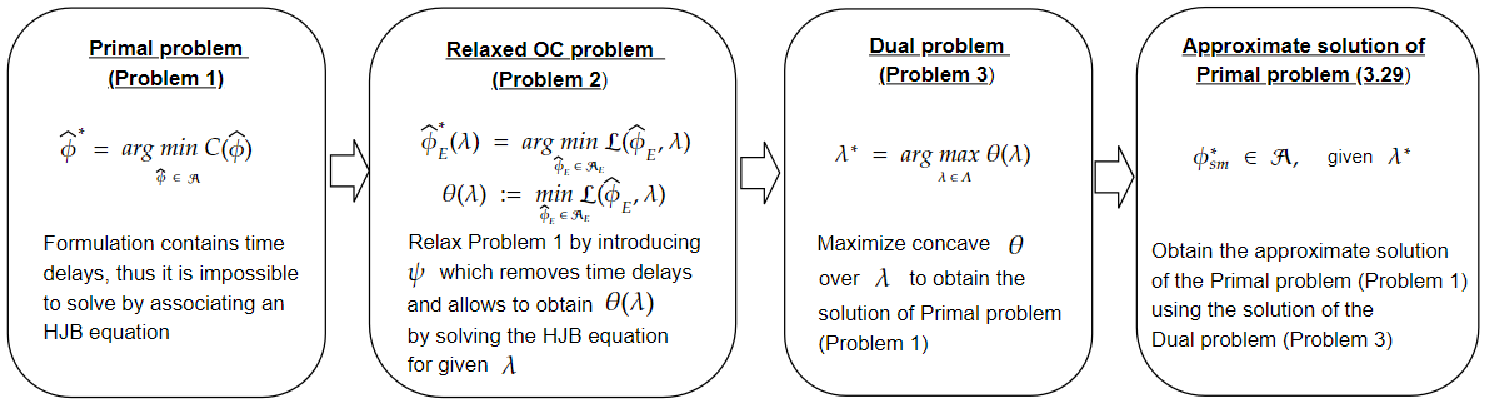}
	\caption{\label{Flowchart} The flowchart showing the connection between Problems~\ref{prob:primal},~\ref{prob:relaxed_oc},~\ref{prob:dual} and the approximate solution of Primal problem stated in Section~\ref{sec:Admissible Controls Smoothing: Computation of Penalized Controls}.}
\end{figure}

In Sections~\ref{sec:Statement of the Optimal Control problem of the Complete Linked System} and~\ref{sec:Continuous-Time Lagrangian Relaxation} we introduce three problems: Problem~\ref{prob:primal} (Primal problem), Problem~\ref{prob:relaxed_oc} (Relaxed OC problem) and Problem~\ref{prob:dual} (Dual problem). Problem~\ref{prob:primal} is the original problem that we want to solve, which aims at finding the optimal controls of operating the power production system minimizing the cost of power production while satisfying the demand and physical constraints. Due to time delays in the formulation of Problem~\ref{prob:primal}, we cannot obtain the solution by solving an associated HJB equation. To address the time delays, we formulate Problem~\ref{prob:relaxed_oc} (Relaxed OC problem) by introducing the virtual controls and relaxing the corresponding constraints by adding the penalization with Lagrangian multipliers of the violation of that constraint into the value function. Problem~\ref{prob:relaxed_oc} does not contain time delays, and thus the optimal control can be obtained by solving the associated HJB equation for given Lagrangian multipliers. We then maximize the dual function of Problem~\ref{prob:relaxed_oc} over Lagrangian multipliers, which constitutes Problem~\ref{prob:dual} (Dual problem). Obtaining optimal Lagrangian multipliers with the corresponding optimal solution allows us to construct the approximate solution of the initial Problem~\ref{prob:primal}, as shown in Section~\ref{sec:Admissible Controls Smoothing: Computation of Penalized Controls}. Figure \ref{Flowchart} shows the connection between the aforementioned problems.

\subsection{Primal Problem}
\label{sec:Statement of the Optimal Control problem of the Complete Linked System}

According to the cost models~\eqref{eq:dams_power_cost_eq} and~\eqref{eq:FFS_power_eq}, 
the instantaneous and total costs associated with the control $\hat{\boldsymbol{\phi}}$ on $[0,T]$ are as follows:
\begin{subequations}
	\label{eq:cost_func_system}
	\begin{align}
		\label{eq:cost_func_instant_system}
		C_I\left(\hat{\boldsymbol{\phi}}(t),\mathbf{X}(t)\right) 
		& = \sum_{i=1}^{N_D}  K_H^{(i)} 
		\underbrace{\left( \bar{\phi}_{Tur}^{(i)}\left(\hat{v}^{(i)}(t)\right)\hat{\phi}_{Tur}^{(i)}(t)
			+\bar{\phi}_S^{(i)}\left(\hat{v}^{(i)}(t)\right) \hat{\phi}_S^{(i)} (t) \right)}_{\hspace{1cm}h^{(i)}\left(\hat{v}^{(i)}(t),\hat{\phi}^{(i)}_{Tur}(t) ,\hat{\phi}^{(i)}_{S}(t)\right)\,:=\hfill}
		+ \sum_{i=1}^{N_F} K_F^{(i)}  \bar{P}_F^{(i)} \hat{\phi}_F^{(i)}(t), \\
		\label{eq:cost_func_total_system}
		C\left(\hat{\boldsymbol{\phi}}\right) &= \int_{0}^{T} C_I\left(\hat{\boldsymbol{\phi}}(t),\mathbf{X}(t)\right) dt. 
	\end{align}
\end{subequations}
The goal is to minimize the total cost $C$ in \eqref{eq:cost_func_total_system}, given 
the dynamics and constraints in Section~\ref{sec:model_summary}.

Then, the OC problem becomes the following.
\begin{problem}[Primal with Time Delays]\label{prob:primal}
	Given the initial data, $\mathbf{X}(0)\in [0,1]^{N_{D}+1}$, and 
	forecasts, $D_E(t)$ and $\left\{I^{(i)}(t)\right\}_{i=1}^{N_D}$, for $0\leq t\leq T$, 
	find 
	\begin{align}
		\hat{\boldsymbol{\phi}}^\ast & = \argmin_{\hat{\boldsymbol{\phi}}}{C\left(\hat{\boldsymbol{\phi}}\right)},
	\end{align}
	where the minimization is taken over all controls that satisfy 
	Constraints~\ref{instantaneous_constraints} and \ref{time_delayed_constraints} 
	with Dynamics~\ref{dynamics}.
\end{problem}
The time delays between the dams in Constraints~\ref{time_delayed_constraints} prevent using the classical dynamic programming algorithm (see, e.g.,~\cite{bertsekas1995dynamic,berkovitz2013optimal}), directly on Problem~\ref{prob:primal}.
One strategy sometimes used in this context is to discretize the time interval, numerically approximate the dynamics, and impose constraints on the control in the discretization points. That strategy leads to a large discrete, constrained optimization problem to which numerical optimization algorithms can be applied directly. With such a strategy, the benefits provided by the structure of an OC problem, in particular the applicability of dynamic programming, are lost.
Other approaches to remove the time delays include state augmentation techniques enlarging the state space. However, the resulting reformulated problem may have complicated a state and control, suffering from the curse of dimensionality.
This work proposes using Lagrangian relaxation over continuous-time constraints, which links the states at different times. Section~\ref{sec:Continuous-Time Lagrangian Relaxation} explains this idea and how it is adapted to the considered problem.

\subsection{Continuous-Time Lagrangian Relaxation and the Dual Problem}\label{sec:Continuous-Time Lagrangian Relaxation}

The continuous-time Lagrangian relaxation is based on the classical Lagrangian relaxation technique applied in constrained optimization problems. This approach allows solving an easier, dual problem 
(Problem~\ref{prob:primal}) by relaxing the constraints using Lagrange multipliers. 
In this work, we adapt this idea to the OC problem to address the delayed dynamics in continuous time. The relaxed formulation enables stating an HJB equation~\cite{fleming2006controlled} and solving the problem using appropriate numerical techniques. 

We let $\mathcal{B}_U:=\cup_{i=1}^{N_D} \mathcal{B}(i)$ denote the set of all upstream dams.
We introduce one virtual control, $\psi^{(j)}(t)$, for each $j\in\mathcal{B}_U$, 
to match the corresponding time-delayed flow.  
The normalization constant is defined as
\begin{align}
	\bar{\psi}^{(j)} &= \underset{\hat{v} \in [0,1]}{\max}  \left( \bar{\phi}_{Tur}^{(j)}(\hat{v}) + \bar{\phi}_S^{(j)}(\hat{v})\right), 
	&& \text{for $j \in \mathcal{B}_U$,}
\end{align}
with the corresponding normalized virtual controls
\begin{align}
	\hat{\psi}^{(j)}(t) &= \frac{\psi^{(j)}(t)}{\bar{\psi}^{(j)}}, 
	&& \text{for $j \in \mathcal{B}_U$,}
\end{align}
The normalized virtual controls are denoted jointly by
$\boldsymbol{\hat\psi}(t):=\left\{ \hat\psi^{(j)}(t): j \in \mathcal{B}_U\right\}$. 
The following two constraints express the goal of matching the time-delayed flows and the normalization of the virtual controls. (The first constraint will be subject to Lagrangian relaxation below.)
\begin{constraints}[Virtual Control Balance]\label{virtual_controls_target}
	For any $j \in \mathcal{B}_U$, 
	\begin{align}
		\bar{\psi}^{(j)} \hat{\psi}^{(j)}(t) &= h^{(j)}(\hat{v}^{(j)}(t-\tau_j),\hat{\phi}^{(j)}_{Tur}(t-\tau_j),\hat{\phi}^{(j)}_{S}(t-\tau_j)), && \text{for $\tau_j<t<T$}. 
	\end{align}
\end{constraints}
The virtual controls are added to the system to replace the time-delayed real water coming from an upstream dam in Dynamics~\ref{dynamics} and Constraints~\ref{time_delayed_constraints} with virtual water at current time, so that the classical dynamic programming algorithm can be used to solve the relaxed OC problem.
\begin{constraints}[Virtual Control Range]\label{virtual_controls_range}
	For any $j \in \mathcal{B}_U$, 
	$0 \leq \hat{\psi}^{(j)}(t) \leq 1$, with $t\in[\tau_j,T]$.
\end{constraints}
After the relaxation of Constraints~\ref{virtual_controls_target} the range constraints only 
guarantee consistency with the maximal possible outflows of the dams, not consistency with the maxima, considering the volume $\hat{v}^{(j)}(t-\tau_j)$. Replacing the virtual controls satisfying Constraints~\ref{virtual_controls_target} in Constraints~\ref{time_delayed_constraints} leads to the following.
\begin{constraints}[Relaxed Control Constraints]\label{relaxed_constraints}
	For any $i=1,2,\dots,N_D$ s.t. $\mathcal{B}(i)$ is nonempty, 
	\begin{subequations}\label{eq:rel_con}
		\begin{align}
			h^{(i)}(1, \bar{\phi}^{(i)}_{Tur}(t),\bar{\phi}^{(i)}_{S}(t))
			& \geq 
			I^{(i)}(t)+\sum_{j\in\mathcal{B}(i)}\bar{\psi}^{(j)} \hat{\psi}^{(j)}(t),  
			&& \text{if $\hat{v}^{(i)} (t)= 1$,}\\
			h^{(i)}(0, \bar{\phi}^{(i)}_{Tur}(t),\bar{\phi}^{(i)}_{S}(t))
			& \leq
			I^{(i)}(t)+\sum_{j\in\mathcal{B}(i)}\bar{\psi}^{(j)} \hat{\psi}^{(j)}(t),  
			&& \text{if $\hat{v}^{(i)}(t) = 0$.}
		\end{align}
	\end{subequations}
\end{constraints}
The extended controls are denoted by $\hat{\boldsymbol{\phi}}_{E}$, where 
$\hat{\boldsymbol{\phi}}_{E}(t) = (\hat{\boldsymbol{\phi}}(t), \hat{\boldsymbol{\psi}}(t)) \in \mathcal{A}_{E}(t, \mathbf{X}(t))$. 
The set of admissible controls, $\mathcal{A}_{E}(t, \mathbf{X}(t))$, is given by Constraints~\ref{instantaneous_constraints}, \ref{virtual_controls_range}, and~\ref{relaxed_constraints}.
We also use $\hat{\boldsymbol{\phi}}_{E}\in\mathcal{A}_{E}$ as shorthand for 
$\hat{\boldsymbol{\phi}}_{E}(t) \in \mathcal{A}_{E}(t, \mathbf{X}(t)), \; \forall t \in [0,T]$. 
Combining the extended controls with Dynamics~\ref{dynamics} leads to a relaxed version of the dynamics and constraints.
\begin{dynamics}[Relaxed Resource Dynamics]\label{relaxed_dynamics}
	The normalized volumes, 
	$\hat{v}^{(i)}$, for $i=1,\dots,N_D$, and battery capacity, $\hat{a}$, 
	satisfy the controlled dynamics
	\begin{align}
		d\hat{v}^{(i)}(t) & = 
		\underbrace{\frac{I^{(i)}(t) - h^{(i)}(\hat{v}^{(i)}(t),\hat{\phi}^{(i)}_{Tur}(t),\hat{\phi}^{(i)}_{S}(t))}{\bar{v}^{(i)}-\underline{v}^{(i)}}}_{\hspace{1cm}g_H^{(i)}\left(t,\hat{v}^{(i)}(t),\hat{\boldsymbol{\phi}}(t) \right)\,:=\hfill} 
		+ \frac{1}{\bar{v}^{(i)}-\underline{v}^{(i)}}\sum_{j\in\mathcal{B}(i)} \bar{\psi}^{(j)} \hat{\psi}^{(j)}(t) 
		dt,\\ 
		d\hat{a}(t)       & = -\frac{\overline{P}_A}{\bar{A}} \hat{\phi}_A(t) dt,
	\end{align}
	in $[0,T]$, with all initial conditions, $\hat{v}^{(i)}(0)$ and $\hat{a}(0)$, in $[0,1]$.
\end{dynamics}
To formulate a dual problem associated with Problem~\ref{prob:primal}, 
we relax Constraints~\ref{virtual_controls_target} by introducing bounded time-continuous deterministic Lagrange multipliers, 
$\boldsymbol{\lambda}(t)=\left\{\lambda_j (t): j \in \mathcal{B}_U\right\}$, 
where each $\lambda_j: [\tau_{j}, T] \mapsto \R$ is associated with one of the constraints in Constraints~\ref{virtual_controls_target}.
For each $j\in \mathcal{B}_U$, we add the following integral to the cost functional in~\eqref{eq:cost_func_total_system}:
\begin{multline}
	\int_{\tau_{j}}^{T} \lambda_j (t) 
	\left( \bar{\psi}^{(j)} \hat{\psi}^{(j)}(t) - 
	h^{(j)}\left(\hat{v}^{(j)}(t-\tau_j),\hat{\phi}^{(j)}_{Tur}(t-\tau_j),\hat{\phi}^{(j)}_{S}(t-\tau_j)\right) \right)dt
	\\ =
	\int_{\tau_{j}}^{T} \lambda_j (t) \bar{\psi}^{(j)} \hat{\psi}^{(j)}(t)\, dt 
	- 
	\int_{0}^{T-\tau_{j}} \lambda_j (t+\tau_{j})
	h^{(j)}\left(\hat{v}^{(j)}(t),\hat{\phi}^{(j)}_{Tur}(t),\hat{\phi}^{(j)}_{S}(t)\right)dt
\end{multline}
to obtain the Lagrangian
\begin{multline}
	\mathcal{L}(\hat{\boldsymbol{\phi}}_{E}, \boldsymbol{\lambda})
	= C\left(\hat{\boldsymbol{\phi}}\right) 
	+ \sum_{j \in \mathcal{B}_U} \int_{\tau_{j}}^{T} \lambda_j (t) \bar{\psi}^{(j)} \hat{\psi}^{(j)}(t)\, dt \\
	- \sum_{j \in \mathcal{B}_U} \int_{0}^{T-\tau_{j}} \lambda_j (t+\tau_{j}) 
	h^{(j)}\left(\hat{v}^{(j)}(t),\hat{\phi}^{(j)}_{Tur}(t),\hat{\phi}^{(j)}_{S}(t)\right)dt.
\end{multline} 
The set of allowed Lagrange multipliers is denoted as $\Lambda$.

\begin{problem}[Relaxed OC Problem]\label{prob:relaxed_oc}
	Given the initial data, $\mathbf{X}(0)\in [0,1]^{N_{D}+1}$, 
	Lagrange multipliers, $\boldsymbol{\lambda}(t)$, and forecasts, $D_E(t)$ and $\left\{I^{(i)}(t)\right\}_{i=1}^{N_D}$, for $0\leq t\leq T$, find 
	\begin{align}
		\hat{\boldsymbol{\phi}}_{E}^\ast{(\cdot,\boldsymbol{\lambda})} & = 
		\underset{\hat{\boldsymbol{\phi}}_{E}\in\mathcal{A}_{E}}{\argmin} \; \mathcal{L}\left(\hat{\boldsymbol{\phi}}_{E}, \boldsymbol{\lambda}\right),
	\end{align}
	where $\hat{\boldsymbol{\phi}}_{E}\in\mathcal{A}_{E}$ indicates that the minimization is taken over all controls that satisfy
	Constraints~\ref{instantaneous_constraints}, \ref{virtual_controls_range}, and \ref{relaxed_constraints} 
	with Dynamics~\ref{relaxed_dynamics}.
\end{problem}
As a consequence of the Lagrange relaxation, the extended controls are no longer coupled in time for a fixed multiplier $\boldsymbol{\lambda}$. Therefore, the dynamic programming principle applies 
to Problem~\ref{prob:relaxed_oc}.
Finally, the dual problem associated with Problem~\ref{prob:primal} is as follows. 
\begin{problem}[Dual Problem]\label{prob:dual}
	Find
	\begin{align}
		\boldsymbol{\lambda}^\ast & = \underset{\boldsymbol{\lambda\in\Lambda} }{\arg\max} \;\theta\left(\boldsymbol{\lambda}\right),
	\end{align}
	where the dual function, $\theta\left(\boldsymbol{\lambda}\right)$, is given by 
	\begin{align}\label{eq:dual function (markovian problem)}
		\theta\left(\boldsymbol{\lambda}\right) &=  
		\mathcal{L}\left(\hat{\boldsymbol{\phi}}_{E}^\ast{(\cdot,\boldsymbol{\lambda})}, \boldsymbol{\lambda}\right),
	\end{align}
	and $\hat{\boldsymbol{\phi}}_{E}^\ast{(\cdot,\boldsymbol{\lambda})}$ solves Problem~\ref{prob:relaxed_oc}, given $\boldsymbol{\lambda}$.
\end{problem}

The structure of the Lagrange multipliers, $\boldsymbol{\lambda}$, or equivalently the set $\Lambda$, depends on the considered power system type. We only introduce Lagrange multipliers in systems with connected dams. In the special case considered in this paper, of deterministic dynamics for the connected dams $\boldsymbol{\lambda}$ are deterministic Lagrange multipliers. In this case, we require only basis functions 
in time to approximate $\boldsymbol{\lambda}$. We provide more details in Section~\ref{sec:Optimization of the Dual Problem}. 
In the more general case of a system of connected dams and stochastic dynamics (e.g., wind or solar), $\boldsymbol{\lambda}$ is the stochastic process that may depend on the whole path of the wind or solar power. For instance, $\boldsymbol{\lambda}$ can be considered the price associated with the hydropower source, given the wind/solar state. The water is more valuable in times of low wind than when the wind is strong. Then, to approximate $\boldsymbol{\lambda}$, we require basis functions in time and w.r.t. realizations ($\omega$). 
The case of stochastic dynamics is left as the subject of future work.

\paragraph{The HJB equation related to the relaxed OC problem.}
An approximation, $\bar{\theta}(\boldsymbol{\lambda})$, of the solution $\theta(\boldsymbol{\lambda})$ to Problem~\ref{prob:relaxed_oc} is obtained by numerically approximating the associated first-order HJB 
initial value problem:
\begin{align}\label{eq:HJB_PDE}
	\begin{cases}
		\frac{\partial u}{\partial t}+H_R\left(t, \hat{\mathbf{v}},\hat{a}, \boldsymbol{\lambda}, Du\right)=0, 
		& t \in [0,T],\: \left[\hat{\mathbf{v}},\hat{a}\right]\in [0,1]^{N_D+1}
		\\
		u(T,\cdot,\cdot; \boldsymbol{\lambda})=0 
	\end{cases}
\end{align}
where the Hamiltonian, $H_R$, is given by
\begin{align}\label{eq:Hamiltonian}
	H_R\left(t, \hat{\mathbf{v}},\hat{a}, \boldsymbol{\lambda}, Du\right) &=
	\underset{{ \hat{\boldsymbol{y}}_{E} } \in \mathcal{A}_{E}(t, \mathbf{X}(t))}{\min} 
	F_R\left({ \hat{\boldsymbol{y}}_{E} };t, \hat{\mathbf{v}},\hat{a}, \boldsymbol{\lambda}, Du\right).
\end{align}
Here, ${\hat{\boldsymbol{y}}_{E} = (\hat{\boldsymbol{y}}, \hat{\boldsymbol{z}})}$  represents an optimization variable with the obvious connections to the control $\hat{\boldsymbol{\phi}}_{E}(t) = (\hat{\boldsymbol{\phi}}(t), \hat{\boldsymbol{\psi}}(t))$  at time $t$ and the same goes for its components; for example, , 
$\hat{y}^{(j)}_{S}$ corresponds to $\hat{\phi}^{(j)}_{S}(t)$.
We recall that ${\hat{\boldsymbol{y}}_{E} } \in \mathcal{A}_{E}(t, \mathbf{X}(t))$ means that it is subject to Constraints~\ref{instantaneous_constraints}, \ref{virtual_controls_range}, 
and \ref{relaxed_constraints} at time $t\in[0,T]$.
With this notation, the objective function is
\begin{multline}
	\label{eq:objective_fun}
	F_R\left({\hat{\boldsymbol{y}}_{E} };t, \hat{\mathbf{v}},\hat{a}, \boldsymbol{\lambda}, Du\right)
	= 
	C_I\left({\hat{\boldsymbol{y}} },\mathbf{X}(t)\right) 
	- \frac{\overline{P}_A}{\bar{A}} {\hat{y}_A} \frac{\partial u}{\partial \hat{a}}(t)
	+ \sum_{i=1}^{N_D} g_H^{(i)}\left(t,\hat{v}^{(i)},{\hat{\boldsymbol{y}} }\right)
	\frac{\partial u}{\partial \hat{v}^{(i)} }(t)\\
	+ \sum_{i=1}^{N_D} \left(\sum_{j\in\mathcal{B}(i)} \frac{\bar{\psi}^{(j)} {\hat{z}^{(j)}}} {\bar{v}^{(i)}-\underline{v}^{(i)}}\right)  \frac{\partial u}{\partial \hat{v}^{(i)} }(t)
	+ \sum_{j\in\mathcal{B}_U} \lambda_j(t) \mathds{1}_{[\tau_{j},T]}(t) \bar{\psi}^{(j)} {\hat{z}^{(j)}}\\
	- \sum_{j\in\mathcal{B}_U} \lambda_j(t+\tau_{j}) \mathds{1}_{[0,T-\tau_{j}]}(t)
	h^{(j)}\left(\hat{v}^{(j)},{\hat{y}^{(j)}_{Tur} }, {\hat{y}^{(j)}_{S} }\right).
\end{multline}
In this function, $\mathds{1}_Y(\cdot)$ denotes the indicator function with support on $Y$, and $Du$ denotes the vector of first-order partial derivatives of $u$ w.r.t. the state-space variables.
Boundary conditions do not need to be imposed due to the characteristics of the deterministic first-order partial differential equation (PDE) in~\eqref{eq:HJB_PDE}. 
Section~\ref{sec: Monotone Scheme for Solving the HJB} discusses the numerical treatment of~\eqref{eq:HJB_PDE}.

\subsubsection{Finite-Dimensional Approximation of the Dual Problem}
\label{par:finite-dim_approx}
The Lagrange multipliers $\boldsymbol{\lambda}(t)=\{\lambda_j (t),\: j \in \mathcal{B}_U\}$ in Problem~\ref{prob:dual} are generally piecewise continuous functions, and Problem~\ref{prob:dual} is an infinite-dimensional optimization problem aiming to satisfy Constraints~\ref{virtual_controls_target} pointwise in time. 
We approximate the infinite-dimensional problem by satisfying 
each finite number of subintervals of $[0, T]$. That is, 
the integral of the pointwise violation of Constraints~\ref{virtual_controls_target}, 
$\bar{\psi}^{(j)} \hat{\psi}^{(j)}(t) - h^{(j)}(\hat{v}^{(j)}(t-\tau_j),\hat{\phi}^{(j)}_{Tur}(t-\tau_j),\hat{\phi}^{(j)}_{S}(t-\tau_j))$, 
over each such subinterval should equal zero. 
This approach corresponds to an approximation of the Lagrange multipliers using piecewise constant functions. 
More precisely, for each
$j \in \mathcal{B}_U$, we define the grid $t_k^{(j)} := \tau_j + k\Delta t_{\lambda_j}$, $k \in \{0,1,...,m_j\}$, 
corresponding to the uniform time step size $\Delta t_{\lambda_j} = (T- \tau_j)/m_j$, and approximate
\begin{align}\label{eq: lambda_approx_first}
	\lambda_j(t) \approx 
	\hat{\lambda}_j(t) 
	& := 
	\sum_{k=0}^{m_j-1} \hat{\alpha}^{(j)}_k \mathbbm{1}_{[t_k^{(j)},t_{k+1}^{(j)})}(t), &&  t \in [\tau_j, T], \: j \in \mathcal{B}_U,
\end{align}
where each $\hat{\pmb{\alpha}}^{(j)}=(\hat{\alpha}^{(j)}_1,\dots,\hat{\alpha}^{(j)}_{m_j})$ is a real vector.
For a given time discretization, we define
\begin{align}
	\label{eq:finite_space_lambda}
	\hat{\Lambda} & = 
	\left\{\boldsymbol{\hat{\lambda}}:
	\text{with $\hat{\lambda}_j$ defined in~\eqref{eq: lambda_approx_first} and $\hat{\pmb{\alpha}}^{(j)}\in\R^{m_j}$, $j\in\mathcal{B}_U$}
	\right\}.
\end{align}
In this approximation, small dual gaps can be obtained even when the time discretizations for 
$\boldsymbol{\hat{\lambda}}$ are much coarser than the time discretizations related to the 
dynamics and controls (see Section~\ref{sec:num_experiments}).

For a simple function $\boldsymbol{\hat{\lambda}}\in\hat{\Lambda}$, the Lagrangian becomes
\begin{align}
	\mathcal{L}(\hat{\boldsymbol{\phi}}_{E}, \boldsymbol{\hat{\lambda}}) = C\left(\hat{\boldsymbol{\phi}}\right)
	+ \sum_{j\in\mathcal{B}_U}\sum_{k=0}^{m_j-1} \hat{\alpha}^{(j)}_{k} \Big(\xi_{\hat{\lambda}_j}\theta\Big)_k,
\end{align}
with
\begin{align}
	\nonumber
	\Big(\xi_{\hat{\lambda}_j}\theta\Big)_k & =
	\int_{t^{(j)}_k}^{t_{k+1}^{(j)}}
	h^{(j)}\left(\hat{v}^{(j)}(t-\tau_j),\hat{\phi}^{(j)}_{Tur}(t-\tau_j),\hat{\phi}^{(j)}_{S}(t-\tau_j)\right) dt \\
	& =
	\label{eq:subgradient}
	\bar{\psi}^{(j)} \int_{t^{(j)}_k}^{t_{k+1}^{(j)}} \hat{\psi}^{(j)}(t) dt
	- \int_{t^{(j)}_k-\tau_j}^{t_{k+1}^{(j)}-\tau_j}
	h^{(j)}\left(\hat{v}^{(j)}(t),\hat{\phi}^{(j)}_{Tur}(t),\hat{\phi}^{(j)}_{S}(t)\right) dt,
\end{align}
for $j\in\mathcal{B}_U$ and $0 \le k\le m_j-1$. 
The subgradient of $\theta(\cdot)$, in~\eqref{eq:dual function (markovian problem)}, at  $\hat{\boldsymbol{\lambda}}$ is  
$\pmb{\xi}_{\hat{\pmb{\lambda}}}\theta: =\{\xi_{\hat{\lambda}_j}\theta, \: j \in \mathcal{B}_U \}$.

From the optimization viewpoint, choosing $\Lambda=\hat{\Lambda}$, the numerical approximation to Problem~\ref{prob:dual}, discussed in Section~\ref{sec:SOC_num_approach} is concave but generally not smooth. Therefore, we use subgradient-based methods to solve Problem~\ref{prob:dual}. 

\begin{remark}[On inequality constraints]
	If having more water in a dam is better, then the equality constraints can be replaced with inequality constraints and signed Lagrange multipliers. This a priori information can be beneficial for the numerical methods. However, determining this information is challenging because completely filling a dam may lead to spillage, decreasing the power of the dam due to the raised downstream water level.
\end{remark}

\begin{remark} [On the interpretation of optimal Lagrange multipliers]
	The optimal Lagrange multipliers $\boldsymbol{\lambda}^\ast$ have a meaningful interpretation. Suppose that, for some $j \in \mathcal{B}_U$, we perturb Constraints~\ref{virtual_controls_target} at one interval $[t_k^{(j)},t_{k+1}^{(j)})$ of the time grid defined in 
	Section~\ref{par:finite-dim_approx} by an amount $\epsilon$ on average, so that
	\begin{gather}\label{cont2}
		\int_{t_k^{(j)}}^{t_{k+1}^{(j)}}\bar{\psi}^{(j)} \hat{\psi}^{(j)}(t) -
		h^{(j)}\left(\hat{v}^{(j)}(t-\tau_j),\hat{\phi}^{(j)}_{Tur}(t-\tau_j),\hat{\phi}^{(j)}_{S}(t-\tau_j)\right) dt = \epsilon ~~~~~
		\text{for $t_k^{(j)}\leq t<t_{k+1}^{(j)}$}. 
	\end{gather}
	Through the dependency of the perturbed OC $\hat{\boldsymbol{\phi}}_E^*(\epsilon)$ on $\epsilon$, the optimal value of the objective function~\eqref{eq:cost_func_total_system} of the primal problem has a parametric dependence on $\epsilon$. We call the resulting optimum 
	$\zeta(\epsilon):=C\left(\hat{\boldsymbol{\phi}}_E^*(\epsilon)\right)$.
	If we differentiate w.r.t. $\epsilon$, at $\epsilon=0$, we obtain (see Section 5.6.3 in \cite{boyd2004convex}, or Section 6.2 in \cite{bazaraa2013nonlinear}) the following:
	\begin{gather}\label{eq:tmp_marginpar}
		\frac{\partial \zeta}{\partial \epsilon}\left(0\right) = - \hat{\lambda}_j^*(t_k^{(j)}).
	\end{gather}
	That is, the negative optimal Lagrange multiplier is the derivative of the optimal objective function $\zeta(\epsilon)$
	w.r.t. the perturbation of the constraint on the virtual controls, indicating how much the cost can be reduced if more water arrives at a downstream dam on the time interval $[t_k^{(j)},t_{k+1}^{(j)})$. With $\epsilon = 0$, the water $h^{(j)}\left(\hat{v}^{(j)}(t-\tau_j),\hat{\phi}^{(j)}_{Tur}(t-\tau_j),\hat{\phi}^{(j)}_{S}(t-\tau_j)\right)$
	sent at time $t-\tau_j$ from an upstream dam must match the water received $\bar{\psi}^{(j)} \hat{\psi}^{(j)}(t)$ at time $t$ by a downstream dam precisely. In economics, $\hat{\lambda}_j^*(\epsilon=0)$ is called the shadow price or marginal cost of the resource in question, which is water in this case. After relaxing the constraint, the following analogy applies. Suppose that instead of water flowing down the river, the upstream dam sells water at time $t-\tau_j$ to a third party and that the downstream dam buys water from the same party at time $t$. The term $-\lambda_j^\ast(t+\tau_{j}) \mathds{1}_{[0, T-\tau_{j}]}(t) h^{(j)}\left(\hat{v}^{(j)}(t),\hat{\phi}^{*(j)}_{Tur}(t),\hat{\phi}^{*(j)}_{S}(t)\right)$
	in the objective function~\eqref{eq:objective_fun}, shifted by $-\tau_j$, can be interpreted as the profit gained by selling the water by the upstream dam at time $t-\tau_j$ for the optimal price $\lambda_j^*(t)$. In addition, the term $\lambda_j^\ast(t) \mathds{1}_{[\tau_{j},T]}(t)\bar{\psi}^{(j)} \hat{\psi}^{(j)}(t)$ represents the cost of buying water by the downstream dam at time $t$ (sold by the upstream dam at time $t-\tau_j$) at the optimal price $\lambda_j^*(t)$. The optimal price $\lambda_j^*(t)$ ensures that the volume of water sold by the upstream dam and bought by the downstream dam match. As demonstrated in the next section, this idea is key for constructing an admissible solution to the primal problem.  
\end{remark}

\subsection{Construction of Nearly Optimal Admissible Controls for Problem~\ref{prob:primal}}\label{sec:Construction of the admissible Controls for the Primal Problem}

In order to construct nearly optimal, admissible, controls for Problem~\ref{prob:primal} based on the solution to its relaxed, dual, Problem~\ref{prob:dual}, 
with corresponding control-state paths $\left(\boldsymbol{\phi}_E^\ast,\mathbf{X}_E^\ast\right)$,
the idea is to consider the dual function $\theta$ evaluated at the optimal Lagrange multipliers, $\boldsymbol{\lambda}^\ast(t)$, and  enforce Constraints~\ref{virtual_controls_target} by substituting the virtual controls
$\psi^{(j)}(t)$ with 
$h^{(j)}\left(\hat{v}^{(j)}(t-\tau_j),\hat{\phi}^{(j)}_{Tur}(t-\tau_j),\hat{\phi}^{(j)}_{S}(t-\tau_j)\right)$
for $ \tau_j<t<T$, $j \in  \mathcal{B}_U$. 
This substitution is justified by the interpretation of the optimal Lagrange multipliers given earlier. The optimal price $\pmb{\lambda}^*(t)$ guarantees that the amount of water sold by the upstream dam and bought by the downstream dam nearly match. At time $t$, we must only decide how much water to sell for the upstream dam, and the optimal price $\pmb{\lambda}^*(t)$ ensures that the optimal decision for the downstream dam is to buy the same amount of water at time $t+\tau_j$.

We compute the admissible control-state paths, $\left(\boldsymbol{\phi}_\mathrm{ad}^\ast,\mathbf{X}_\mathrm{ad}^\ast\right)$, 
corresponding to $\boldsymbol{\lambda}^\ast$, forward in time following Dynamics~\ref{dynamics}. The state vector $\mathbf{X}_\mathrm{ad}^\ast = \left({v}^{^\ast(1)},\dots,{v}^{^\ast(N_D)},{a}^\ast\right)$ and normalized controls $\boldsymbol{\phi}_\mathrm{ad}^\ast$ in the flux
are obtained at time $t$ as minimizers of an objective function corresponding to a modified Hamiltonian using the derivatives of the value function 
$\{\frac{\partial u}{\partial \hat{v}^{(i)} }\}_{i=1}^{N^D}$ and $ \frac{\partial u}{\partial \hat{a}}$ obtained when solving the dual problem.
At time $t\in[0, T]$, this objective function 
is viewed as a function of the control at time $t$, depending on data at time $t$, namely $t$, $\mathbf{X}_\mathrm{ad}^\ast(t)$, $Du\left(t,\mathbf{X}_\mathrm{ad}^\ast(t)\right)$, and $\boldsymbol{\lambda}^\ast(t)$. 
Furthermore, the data defining the function include the Lagrange multiplier 
$\boldsymbol{\lambda}^\ast$ at later times $\mathbf{t^{+}} := \cup_{j\in\mathcal{B}_U} \{t+\tau_{j}\} \cap [t,T]$ and $\boldsymbol{\phi}_\mathrm{ad}^\ast$ and $\mathbf{X}_\mathrm{ad}^\ast$ at earlier times $\mathbf{t^{-}} := \cup_{j\in\mathcal{B}_U} \{t-\tau_{j}\} \cap [0,t].$ In the objective function for the control at time $t$,
the $\boldsymbol{\phi}_{\mathrm{ad}}^{\ast(j)}(t-\tau_{j})$ are not controls but inputs determined at earlier times, defined to take zero values in the corresponding time intervals $[0, \tau_j]$, as described after equation~\eqref{eq: dams_connected_constraints}. 
The objective function is
\begin{multline}\label{eq:objective_modified}
	F_\mathrm{ad}
	\left(\boldsymbol{\phi}_\mathrm{ad}^\ast(t) \, ; \, t,\mathbf{X}_\mathrm{ad}^\ast(t), 
	Du\left(t,\mathbf{X}_\mathrm{ad}^\ast(t)\right), \left\{\boldsymbol{\lambda}^\ast(s)\right\}_{s\in t\cup\mathbf{t^{+}}},
	\left\{\boldsymbol{\phi}_\mathrm{ad}^\ast(s),\mathbf{X}_\mathrm{ad}^\ast(s)\right\}_{s\in\mathbf{t^{-}}}
	\right)
	= 
	C_I\left(\boldsymbol{\phi}_\mathrm{ad}^\ast(t)\right)\\ 
	- \frac{\overline{P}_A}{\bar{A}} \boldsymbol{\phi}_{\mathrm{ad},A}^\ast(t) \frac{\partial u}{\partial \hat{a}}(t,\mathbf{X}_\mathrm{ad}^\ast(t))
	+ \sum_{i=1}^{N_D} g_H^{(i)}\left(t,{v}^{\ast(i)},\boldsymbol{\phi}_\mathrm{ad}^\ast(t)\right)
	\frac{\partial u}{\partial \hat{v}^{(i)} }(t,\mathbf{X}_\mathrm{ad}^\ast(t))\\
	+ \sum_{i=1}^{N_D} 
	\left(\sum_{j\in\mathcal{B}(i)}  \frac{h^{(j)}\left(v^{\ast(j)}(t-\tau_j),\phi^{\ast(j)}_{\mathrm{ad},Tur}(t-\tau_j),\phi^{\ast(j)}_{\mathrm{ad},S}(t-\tau_j)\right)
	}{\bar{v}^{(i)} -  \underline{v}^{(i)}}\right) 
	\frac{\partial u}{\partial \hat{v}^{(i)} }(t,\mathbf{X}_\mathrm{ad}^\ast(t))\\
	+ \sum_{j\in\mathcal{B}_U} \lambda_j^\ast(t) \mathds{1}_{[\tau_{j},T]}(t) 
	h^{(j)}\left(v^{\ast(j)}(t-\tau_j),\phi^{\ast(j)}_{\mathrm{ad},Tur}(t-\tau_j),\phi^{\ast(j)}_{\mathrm{ad},S}(t-\tau_j)\right)\\
	- \sum_{j\in\mathcal{B}_U} \lambda_j^\ast(t+\tau_{j}) \mathds{1}_{[0,T-\tau_{j}]}(t) 
	h^{(j)}\left(v^{\ast(j)}(t),\phi^{\ast(j)}_{\mathrm{ad},Tur}(t),\phi^{\ast(j)}_{\mathrm{ad},S}(t)\right).
\end{multline}
In addition, the modified Hamiltonian is
\begin{multline}\label{eq:Hamiltonian_primal}
	{H}_\mathrm{ad}^\ast\left(t,\mathbf{X}_\mathrm{ad}^\ast(t)\, ; \,
	Du\left(t,\mathbf{X}_\mathrm{ad}^\ast(t)\right), \left\{\boldsymbol{\lambda}^\ast(s)\right\}_{s\in t\cup\mathbf{t^{+}}},
	\left\{\boldsymbol{\phi}_\mathrm{ad}^\ast(s),\mathbf{X}_\mathrm{ad}^\ast(s)\right\}_{s\in\mathbf{t^{-}}}
	\right)
	\\ = 
	\underset{\boldsymbol{\phi}_\mathrm{ad}^\ast(t) \in \mathcal{A}(t, \mathbf{X}_\mathrm{ad}^\ast(t))}{\min} 
	F_\mathrm{ad}\left(\boldsymbol{\phi}_\mathrm{ad}^\ast(t)\, ; \,t,\mathbf{X}_\mathrm{ad}^\ast(t), 
	Du\left(t,\mathbf{X}_\mathrm{ad}^\ast(t)\right), \left\{\boldsymbol{\lambda}^\ast(s)\right\}_{s\in t\cup\mathbf{t^{+}}},
	\left\{\boldsymbol{\phi}_\mathrm{ad}^\ast(s),\mathbf{X}_\mathrm{ad}^\ast(s)\right\}_{s\in\mathbf{t^{-}}}
	\right),
\end{multline}
where the space of admissible controls, $\mathcal{A}(t, \mathbf{X}_\mathrm{ad}^\ast(t))$, is characterized by Constraints~\ref{instantaneous_constraints} and Constraints~\ref{time_delayed_constraints}.
Using the obtained admissible controls in the primal cost function~\eqref{eq:cost_func_total_system}, $C\left(\boldsymbol{\phi}_\mathrm{ad}^\ast\right)$, for $\boldsymbol{\phi}_\mathrm{ad}^\ast\in \mathcal{A}$, given $\boldsymbol{\lambda}^\ast$, 
results in a nearly optimal primal solution (i.e., an upper bound for the exact primal solution). 

\subsection{Smoothing of Nearly Optimal Controls}\label{sec:Admissible Controls Smoothing: Computation of Penalized Controls}
The nearly optimal admissible controls for the primal problem constructed in Section~\ref{sec:Construction of the admissible Controls for the Primal Problem} may have time variations that are too fast for practical implementation by the operator. To motivate the approach to smoothing the controls, 
we consider a penalized objective function
\begin{equation} \label{eq: penalized_value_function}
	C(\pmb{\phi}) + \pmb{\beta}^T \pmb{V}(\pmb{\phi};\pmb{\lambda},\pmb{\beta}),
\end{equation}
where the control $\pmb{\phi}$ depends on $\pmb{\lambda}$ and $\pmb{\beta}$. 
In addition, $\boldsymbol{\beta}:=\left(\{\beta^{(i)}_{Tur}, \beta^{(i)}_S\}_{i=1}^{N_D}, \{\beta^{(i)}_F\}_{i=1}^{N_F}, \beta_A \right)$ is the vector of nonnegative penalizing coefficients and
\begin{equation} \label{eq: variation_term}
	\pmb{V}(\pmb{\phi};\pmb{\lambda},\pmb{\beta}) := \Big( \Bigl\{ \int_0^T (\partial_t \phi^{(i)}_{Tur})^2dt,\int_0^T (\partial_t \phi^{(i)}_S)^2dt \Bigr\}_{i=1}^{N_D}, \Bigl\{\int_0^T (\partial_t \phi^{(i)}_F)^2dt \Bigr\}_{i=1}^{N_F},\int_0^T (\partial_t \phi_A)^2dt\Big)
\end{equation}
is the vector whose elements are the squared L$^2$-norms of the time derivatives of the controls.

Because of the derivatives of the controls, we cannot directly associate an HJB equation to the minimization of~\eqref{eq: penalized_value_function}. However, the form of this penalized objective function provides intuition on how to construct a smoothed admissible path. We compute the control-state paths, $\left(\boldsymbol{\phi}_\mathrm{sm}^\ast,\mathbf{X}_\mathrm{sm}^\ast\right)$, 
corresponding to $\boldsymbol{\lambda}^\ast$, forward in time again following Dynamics~\ref{dynamics}, for $\mathbf{X}_\mathrm{sm}^\ast = \left({v}^{^\ast(1)},\dots,{v}^{^\ast(N_D)},{a}^\ast\right)$, where the normalized controls $\boldsymbol{\phi}_\mathrm{sm}^\ast$ in the flux
are obtained at time $t$ as minimizers of an objective function corresponding to another modified Hamiltonian. 
We approximate the derivatives of the controls in~\eqref{eq: penalized_value_function} using finite differences,
(e.g.,
$\partial_t \phi (t)\approx \frac{\phi(t)-\phi(t-\delta t)}{\delta t}=: D_{\delta t}\left[\phi\right](t)$), 
where $\delta t$ is related to the numerical approximation scheme.\footnote{In practice, if the numerical approximation uses nonuniform time discretizations, several different values of $\delta t$ can be used, but we omit this complication here.}
This approach leads to a penalized version of~\eqref{eq:objective_modified} in
\begin{multline}
	F_\mathrm{sm}\left(\boldsymbol{\phi}_\mathrm{sm}^\ast(t) \, ; \, t,\mathbf{X}_\mathrm{sm}^\ast(t), \pmb{\beta}, 
	Du\left(t,\mathbf{X}_\mathrm{sm}^\ast(t)\right), \left\{\boldsymbol{\lambda}^\ast(s)\right\}_{s\in t\cup\mathbf{t^{+}}},
	\left\{\boldsymbol{\phi}_\mathrm{sm}^\ast(s),\mathbf{X}_\mathrm{sm}^\ast(s)\right\}_{s\in\mathbf{t^{-}}\cup\{t-\delta t\}}
	\right)
	\\ = 
	F_\mathrm{ad}\left(\boldsymbol{\phi}_\mathrm{sm}^\ast(t) \, ; \, t,\mathbf{X}_\mathrm{sm}^\ast(t), 
	Du\left(t,\mathbf{X}_\mathrm{sm}^\ast(t)\right), \left\{\boldsymbol{\lambda}^\ast(s)\right\}_{s\in t\cup\mathbf{t^{+}}},
	\left\{\boldsymbol{\phi}_\mathrm{sm}^\ast(s),\mathbf{X}_\mathrm{sm}^\ast(s)\right\}_{s\in\mathbf{t^{-}}}
	\right)
	\\ + 
	\sum_{i=1}^{N_D} 
	\beta_{Tur}^{(i)}
	D^2_{\delta t}\left[\overline{\phi}^{(i)}_{Tur}(\hat{v}^{(i)}(\cdot))\phi^{\ast(i)}_{\mathrm{sm},Tur}\right](t) +
	\sum_{i=1}^{N_D} 
	\beta_{S}^{(i)}
	D^2_{\delta t}\left[\overline{\phi}^{(i)}_{S}(\hat{v}^{(i)}(\cdot))\phi^{\ast(i)}_{\mathrm{sm},S}\right](t) \\
	+ \sum_{i=1}^{N_F} 
	\beta_F^{(i)} 
	D^2_{\delta t}\left[\overline{P}^{(i)}_F\phi^{\ast(i)}_{\mathrm{sm},F}\right](t)
	+ \beta_A
	D^2_{\delta t}\left[\overline{\phi}_A\phi^{\ast}_{\mathrm{sm},A}\right](t),
\end{multline}
corresponding to the penalized Hamiltonian 
\begin{align}\label{eq:Hamiltonian_smooth}
	{H}_\mathrm{sm}^\ast\left(t,\mathbf{X}_\mathrm{sm}^\ast(t)\, ; \,\pmb{\beta}, 
	Du\left(t,\mathbf{X}_\mathrm{sm}^\ast(t)\right), \left\{\boldsymbol{\lambda}^\ast(s)\right\}_{s\in t\cup\mathbf{t^{+}}},
	\left\{\boldsymbol{\phi}_\mathrm{sm}^\ast(s),\mathbf{X}_\mathrm{sm}^\ast(s)\right\}_{s\in\mathbf{t^{-}}\cup\{t-\delta t\}}
	\right)
	\\ = \nonumber
	\underset{\boldsymbol{\phi}_\mathrm{sm}^\ast(t) \in \mathcal{A}(t, \mathbf{X}_\mathrm{sm}^\ast(t))}{\min} 
	F_\mathrm{sm}\left(\boldsymbol{\phi}_\mathrm{sm}^\ast(t)\, ; \,t,\mathbf{X}_\mathrm{sm}^\ast(t), \pmb{\beta}, 
	Du\left(t,\mathbf{X}_\mathrm{sm}^\ast(t)\right), \left\{\boldsymbol{\lambda}^\ast(s)\right\}_{s\in t\cup\mathbf{t^{+}}},
	\left\{\boldsymbol{\phi}_\mathrm{sm}^\ast(s),\mathbf{X}_\mathrm{sm}^\ast(s)\right\}_{s\in\mathbf{t^{-}}\cup\{t-\delta t\}}
	\right).
\end{align}
By varying $\pmb{\beta} > 0$, we determine the control path, which is close to the optimal path but smoother i.e., a smoothed solution of the original primal problem, Problem~\ref{prob:primal}). All controls do not need to be penalized--only those that are impractical. Section~\ref{sec:Admissible Controls Smoothing: The Penalization Parameter Tuning Procedure} explains the procedure related to the tuning of the penalizing coefficients $\boldsymbol{\beta}$.

Similarly to Section~\ref{sec:Construction of the admissible Controls for the Primal Problem}, using the smoothed admissible controls in the primal cost function~\eqref{eq:cost_func_total_system} results in a nearly optimal primal solution, that is, an upper bound for the exact primal solution:
\begin{align}\label{eq:cost_smooth}
	C\left(\boldsymbol{\phi}_\mathrm{sm}^\ast;\boldsymbol{\lambda}^\ast,\pmb{\beta}\right) & = C\left(\boldsymbol{\phi}_\mathrm{sm}^\ast\right), 
	&& \text{for $\boldsymbol{\phi}_\mathrm{sm}^\ast\in \mathcal{A}$, given $\boldsymbol{\lambda}^\ast$ and $\pmb{\beta}$.}
\end{align}
\section{Optimal Control and Lagrangian Relaxation:  Numerical Approach and Error Discussion}\label{sec:SOC_num_approach}
This section discusses the numerical approach to solving Problem~\ref{prob:primal}, following the method described in Section~\ref{sec:SOC_formulation}. 
We describe the setting where all of the dependent dams are serially connected along one of the rivers, which is the case for the system modeled in Section~\ref{sec:num_experiments}. Along that river there are $N^D_1$ dams, numbered so that $\mathcal{B}(i)=i-1$, for $i=2,\dots,N^D_1$, and $\mathcal{B}_U=\{1,\dots,N^D_1-1\}$. 
The remaining dams are numbered $N^D_1+1,\dots,N^D$, and all satisfy $\mathcal{B}(i)=\{\}$.

The overall procedure is illustrated in Algorithm~\ref{algo_gen}, and each block is explained in detail in the following sections.
\begin{algorithm}[h!]
	\caption{{Steps of the approach}}
	\label{algo_gen}
	\begin{algorithmic}[1]
		
		\State {Initialize $\widetilde{\pmb{\lambda}}_{(0)} \leftarrow \pmb{\lambda}_0$, $\pmb{\beta} \leftarrow \pmb{\beta}_0$,  $m_j=1,~j \in \mathcal{B}_U$,   $n=1$ (optimization level), and $\text{TOL}$ denotes a prescribed tolerance}
		\While {Relative duality gap (Error II in (\ref{eq:error_bound_dual_gap})) $> \text{TOL}$}
		{\For {$k \gets 1$ to $N_{\text{iter}}$} 
			\State {Approximate $\theta(\widetilde{\pmb{\lambda}}_{(n-1)}^{(k-1)})$ in \eqref{eq:dual function (markovian problem)} with $\widetilde{\theta}(\widetilde{\pmb{\lambda}}_{(n-1)}^{(k-1)})$ by solving the HJB equation \eqref{eq:HJB_PDE} 
				as explained in Section~\ref{sec: Monotone Scheme for Solving the HJB}.} 
			\State {Compute subgradients 	$\boldsymbol{\xi}_{\widetilde{\pmb{\lambda}}_{(n-1)}^{(k-1)}}\theta$}
			\State {$\widetilde{\pmb{\lambda}}_{(n)}^{(k)} \leftarrow \mathcal{G}(\widetilde{\theta}(\widetilde{\pmb{\lambda}}_{(n-1)}^{(k-1)}), \boldsymbol{\xi}_{\widetilde{\pmb{\lambda}}_{(n-1)}^{(k-1)}}\widetilde{\theta})$:}   \quad updating rule of the LMBM (see Section~\ref{sec:Optimization of the Dual Problem} and Appendix~\ref{sec:Formulation of the LMBM Optimization Method})
			\EndFor}
		\State{\textbf{end}}
		\State {Compute $\widetilde{C}(\widetilde{\pmb{\phi}}^{\ast}_\mathrm{sm}; \widetilde{\pmb{\lambda}}_{(n)}^{\ast}, \boldsymbol{\beta})$, an approximation of the admissible smooth cost~\eqref{eq:cost_smooth}
			of Problem~\ref{prob:primal} as explained in Sections~\ref{sec:Admissible Controls Smoothing: Computation of Penalized Controls} and~\ref{sec:Admissible Controls Smoothing: The Penalization Parameter Tuning Procedure}.}
		\State { Evaluate the relative duality gap as in Error (II) in \eqref{eq:error_bound_dual_gap} between $\widetilde{C}(\widetilde{\pmb{\phi}}^{\ast}_\mathrm{sm}; \widetilde{\pmb{\lambda}}_{(n)}^{\ast}, \boldsymbol{\beta})$ and $\theta(\widetilde{\pmb{\lambda}}_{(n)}^{\ast})$}
		\State {$n=n+1$, $m_j = 2^{n-1},~j \in \mathcal{B}_U$ and initialize $\widetilde{\pmb{\lambda}}_{(n)}^{(0)}$ using $ \widetilde{\pmb{\lambda}}_{(n-1)}^{\ast}$ and update $\boldsymbol{\beta}$ as in Section~\ref{sec:Admissible Controls Smoothing: The Penalization Parameter Tuning Procedure} }
		\EndWhile
		\State \textbf{end}
	\end{algorithmic}
\end{algorithm}

\subsection{{Monotone Scheme for Solving the HJB Equation}}
\label{sec: Monotone Scheme for Solving the HJB}
{To solve the HJB equation defined by~\eqref{eq:HJB_PDE}--\eqref{eq:Hamiltonian} and obtain the numerical approximation $\widetilde{\theta}(\widetilde{\pmb{\lambda}}_{(n)}^{(k)})$ of the dual function $\theta(\widetilde{\pmb{\lambda}}_{(n)}^{(k)})$ for a given $\widetilde{\pmb{\lambda}}_{(n)}^{(k)}$, we employ the backward-in-time forward Euler upwind finite-difference scheme ~\cite{sun2015HJBscheme}, which provides a convergent approximation to the solution of the HJB equation under certain conditions. In \ref{sec: Monotone Scheme for Solving the HJB}, we explain the details of the scheme adapted to this context.
	
	Upon solving the HJB equation with a given $\widetilde{\pmb{\lambda}}_{(n)}^{(k)}$, we have $U_{\pmb{i}}^j$ at every point of the grid $(t_j,\pmb{x}_{\pmb{i}}) \in \mathcal{D}$. This approach allows computing the numerical OC-state path $(\widetilde{\pmb{\phi}}^{\ast}_E,\widetilde{\pmb{X}}^{\ast})$ forward in time starting from the initial conditions $\pmb{X}(0)$ and using Dynamics~\ref{relaxed_dynamics}. The controls that appear in Dynamics~\ref{relaxed_dynamics} are obtained by minimizing (\ref{eq:objective_fun}) at a current point of the state path $\widetilde{\pmb{X}}^{\ast}$, which may be a point outside of the grid $\mathcal{D}$. The values of the derivatives of the value function at those points are interpolated linearly from the derivatives of the value function at the nearest grid points. The subgradient $\pmb{\xi}_{\widetilde{\pmb{\lambda}}_{(n)}^{(k)}}\theta$ is computed according to (\ref{eq:subgradient}) over the OC-state path $(\widetilde{\pmb{\phi}}^{\ast}_E,\widetilde{\pmb{X}}^{\ast})$.
	
	\subsection{Dual Problem Optimization}\label{sec:Optimization of the Dual Problem}
	This section explains how we solve the nonsmooth concave dual optimization (Problem~\ref{prob:dual}) after approximating the dual function for a given $\pmb{\lambda}$, as explained in Section~\ref{sec: Monotone Scheme for Solving the HJB}. 
	The optimization is conducted in two stages: (i) construction of an approximation, $\hat{\boldsymbol{\lambda}}(t)$, of $\boldsymbol{\lambda}(t)$ with basis functions in time, as in (\ref{eq: lambda_approx_first}), and (ii) numerical optimization of the coordinates of the approximated version $\hat{\boldsymbol{\lambda}}(t)$ resulting in the final approximation of the Lagrangian multipliers functions, $\widetilde{\boldsymbol{\lambda}}(t)$. 
	
	\subsubsection{Refinement of the Lagrangian Multiplier Functions}\label{sec:Refinement of the Lagrangian multiplier functions}
	
	We construct the approximation as explained in Section~\ref{sec:Continuous-Time Lagrangian Relaxation}.
	Increasing the level of refinements $\{m_j\}_{j \in \mathcal{B}_U}$ improves the accuracy of the approximation of $ \boldsymbol{\lambda}(t)$ by $ \hat{\boldsymbol{\lambda}}(t)$; however, it increases the dimension of the optimization, which is $\sum_{j \in \mathcal{B}_U} m_j$. In this context, we do not select the number of level refinements based on controlling the approximation error of the Lagrangian multiplier functions but rather on achieving the desired duality gap.
	
	\subsubsection{Numerical Optimization Procedure}\label{sec:The numerical Optimization Procedure}
	
	The second stage of the dual problem optimization is solving the corresponding optimization of $\theta(\hat{\pmb{\lambda}}(t) )$ w.r.t. the parameters $(\hat{\pmb{\alpha}}^{(1)},\dots,\hat{\pmb{\alpha}}^{(N^D_1-1)})$ in~\eqref{eq: lambda_approx_first}.
	Two main types of methods for nonsmooth optimization exist \cite{karmitsa2012comparing, bagirov2014introduction}: (i) subgradient methods (see \cite{kappel2000implementation,beck2003mirror,boyd2003subgradient} and Chapter 3 in \cite{nesterov2018lectures}), and (ii) bundle methods \cite{makela1992nonsmooth,makela2002survey}. The subgradient methods use one arbitrary subgradient at each point, whereas the bundle methods approximate the whole subdifferential of the objective function, making them more powerful in terms of convergence speed for medium- and large-scale problems \cite{karmitsa2012comparing}. In this work, we use the LMBM \cite{LMBM_PhD,haarala2004new,haarala2007globally}, which is a hybrid version of the variable metric bundle methods \cite{lukvsan1999globally,vlvcek2001globally} and the limited memory variable metric methods \cite{byrd1994representations}. This algorithm takes the value of the function at the current point and its subgradient as input. We present a summary of the LMBM method formulation and its algorithm in Appendix~\ref{sec:Formulation of the LMBM Optimization Method}.
	
	At each level of refinement $\{m_j\}_{j \in \mathcal{B}_U}$, we run the LMBM for a fixed number of iterations $N_{\text{iter}}$. In practical experiments, significant improvement in $\widetilde{\theta}(\widetilde{\pmb{\lambda}})$ is achieved after the first few iterations. However, further progress is limited by the accuracy of the current approximation of $\pmb{\lambda}(t)$ (see, e.g., Figure~\ref{NE_LamOpt_R2}). Therefore, the value of $N_{\text{iter}}$ is fixed in Algorithm~\ref{algo_gen} and verified to be sufficiently large.
	
	The result of running $N_{\text{iter}}$ iterations of LMBM optimization at the current level of refinement $n$ is an approximated vector $\widetilde{\pmb{\lambda}}^{\ast}_{(n)} = (\widetilde{\pmb{\alpha}}^{\ast(1)},\dots,\widetilde{\pmb{\alpha}}^{\ast(N^D_1-1)})$. The optimal value $\widetilde{\pmb{\lambda}}^{\ast}_{(n)}$ is not necessarily achieved at the last iteration because the objective function to be minimized does not always decrease at every step in subgradient and bundle type methods. Thus, the best value of the objective function achieved and its corresponding minimizer are stored as optimal in such methods.
	
	\subsection{Admissible Controls Smoothing: Penalization Parameter-Tuning Procedure}
	\label{sec:Admissible Controls Smoothing: The Penalization Parameter Tuning Procedure}
	Upon solving the HJB equation, we aim to obtain the corresponding controls for the admissible smoothed path using \eqref{eq:Hamiltonian_smooth} to compute the numerical smoothed primal cost, ${\widetilde{C}(\widetilde{\pmb{\phi}}^{\ast}_{\mathrm{sm}}; \widetilde{\pmb{\lambda}}^{\ast}_{(n)}, \boldsymbol{\beta})}$. To smooth these controls, we must carefully choose the penalization parameters $\boldsymbol{\beta}$ because large values induce a high primal cost and, consequently, a large dual gap. Moreover, small values do not produce the desired smoothing behavior for the controls. 
	
	All penalization coefficients in $\pmb{\beta}$ are tuned separately using the same method. For illustration, we demonstrate how to choose the coefficient $\beta_A$ to penalize the battery control while maintaining the remaining coefficients in $\pmb{\beta}$ equal to zero.  
	
	The coefficient $\beta_A$ in \eqref{eq: penalized_value_function} can be viewed as a Tikhonov regularization parameter and can be chosen using the notion of the L-curve \cite{hansen1992Lcurve}. As $\beta_A$ increases from zero to $\infty$, the variation of the control term $\widetilde{V}_{\beta_A}(\widetilde{\pmb{\phi}}^{\ast}_{\mathrm{sm}};\widetilde{\pmb{\lambda}}^{\ast}_{(n)},\beta_A)$ in the penalized objective function~\eqref{eq: penalized_value_function} decreases, which is desired. However, the primal cost $\widetilde{C}(\widetilde{\pmb{\phi}}^{\ast}_{\mathrm{sm}};\widetilde{\pmb{\lambda}}^{\ast}_{(n)},\beta_A)$ increases, which should be avoided. Figure~\ref{NE_LCurves} illustrates the compromise between these two conflicting objectives. The L-curve divides the space into two regions: a horizontal region with a small $\beta_A$ where the cost is not increasing too much but the variation becomes significantly reduced, and a vertical region where the controls are already effectively flattened and the variation can no longer be reduced. The primal cost is very high because we deviated too far from the optimal solution. The idea is to choose $\beta_A$, which lies at the leftmost part of the horizontal region of the L-curve (or at the corner, if the curve has a sharp L-shape) to achieve the best effect of smoothing the controls without dramatically increasing the primal cost. Mathematically, for a small $\delta$ value, this can be formulated as follows for some user-specified $\epsilon>0$: 
	\begin{gather}
		\underset{\beta_A \in [0,\infty)}{\max} \beta_A, ~~~~
		\text{s.t.}~~  \left|\frac{\widetilde{C}(\widetilde{\pmb{\phi}}^{\ast}_{\mathrm{sm}}; \widetilde{\pmb{\lambda}}^{\ast}_{(n)}, \beta_A+\delta) - \widetilde{C}(\widetilde{\pmb{\phi}}^{\ast}_{\mathrm{sm}}; \widetilde{\pmb{\lambda}}^{\ast}_{(n)}, \beta_A)}{\widetilde{V}(\widetilde{\pmb{\phi}}^{\ast}_{\mathrm{sm}};\widetilde{\pmb{\lambda}}^{\ast}_{(n)},\beta_A+\delta)-\widetilde{V}(\widetilde{\pmb{\phi}}^{\ast}_{\mathrm{sm}};\widetilde{\pmb{\lambda}}^{\ast}_{(n)},\beta_A)}\right| < \epsilon.
	\end{gather}
	
	\subsection{{Error Discussion of the Approach}}\label{sec:Numerical Analysis of our Approach}
	
	This section explains the different numerical errors involved in the proposed approach. The aim is to obtain an upper bound for the exact dual gap between the primal, $C(\hat{\pmb{\phi}}^\ast)$ (\ref{eq:cost_func_total_system}), and the dual, $\theta(\pmb{\lambda}^\ast)$ (\ref{eq:dual function (markovian problem)}), costs corresponding to using the exact optimal continuous controls of Problem~\ref{prob:primal}, 
	$\hat{\pmb{\phi}}^{\ast}$,	 
	and exact optimal Lagrangian multiplier functions 
	$\pmb{\lambda}^{\ast}$. We let 
	$\widetilde{\pmb{\phi}}^{\ast}_{\mathrm{sm}}$ and 
	$\widetilde{\pmb{\lambda}}$ be the numerical solutions of our approach. Then, by the optimality condition, we have 
	\begin{small}
		\begin{equation} \label{eq: Dual cost continuous}
			\text{Exact Dual Gap}(\pmb{\lambda}^\ast,\hat{\pmb{\phi}}^\ast):=  C(\hat{\pmb{\phi}}^{\ast})- \theta (\pmb{\lambda}^\ast) \leq  C(\widetilde{\pmb{\phi}}^{\ast}_{\mathrm{sm}})- \theta (\widetilde{\boldsymbol{\lambda}}).  
		\end{equation}
	\end{small}
	We let $\widetilde{C}(\widetilde{\pmb{\phi}}^{\ast}_{\mathrm{sm}};\widetilde{\boldsymbol{\lambda}},\pmb{\beta})$ and $\widetilde{\theta}(\widetilde{\boldsymbol{\lambda}})$ denote the numerical approximations of the primal admissible smooth and dual cost, respectively, obtained via Algorithm~\ref{algo_gen}. Moreover, $u(t_0,\pmb{X}_0;\pmb{\lambda},\Delta t, \Delta \pmb{X})$ denotes the HJB solution on a given grid with mesh size ($\Delta t, \Delta \pmb{X}$). Then, we have the following relative error decomposition:
	\begin{align}\label{eq:error_bound_dual_gap}
		\frac{\left| C(\widetilde{\pmb{\phi}}^{\ast}_{\mathrm{sm}})- \theta (\widetilde{\boldsymbol{\lambda}}) \right| }{\widetilde{\theta}(\widetilde{\boldsymbol{\lambda}})}
		&	\leq \Bigg( \underset{\text{Error I: Primal error}}{\underbrace{\frac{|C(\widetilde{\pmb{\phi}}^{\ast}_{\mathrm{sm}}) - \widetilde{C}(\widetilde{\pmb{\phi}}^{\ast}_{\mathrm{sm}};\widetilde{\boldsymbol{\lambda}},\pmb{\beta}) |}{\widetilde{\theta}(\widetilde{\boldsymbol{\lambda}})}}}+ \underset{\text{Error II: Numerical Dual gap}}{\underbrace{\frac{|\widetilde{C}(\widetilde{\pmb{\phi}}^{\ast}_{\mathrm{sm}};\widetilde{\boldsymbol{\lambda}},\pmb{\beta})  - \widetilde{\theta}(\widetilde{\boldsymbol{\lambda}})|}{\widetilde{\theta}(\widetilde{\boldsymbol{\lambda}})}}}  \nonumber\\
		&+ \underset{\text{Error III: Dual approximation error}}{\underbrace{\frac{|\widetilde{\theta}(\widetilde{\boldsymbol{\lambda}}) - u(t_0,\pmb{X}_0;\widetilde{\pmb{\lambda}},\Delta t, \Delta \pmb{X}) |}{\widetilde{\theta}(\widetilde{\pmb{\lambda}})}}}
		+ \underset{\text{Error IV:  HJB error}}{\underbrace{\frac{| u(t_0,\pmb{X}_0;\widetilde{\boldsymbol{\lambda}},\Delta t, \Delta \pmb{X}) - \theta(\widetilde{\boldsymbol{\lambda}})|}{\widetilde{\theta}(\widetilde{\boldsymbol{\lambda}})}}} \Bigg)
	\end{align}
	The numerical dual gap (Error II) and dual approximation (Error III) errors are quantities that we can directly measure. The primal (Error I) and HJB (Error IV) errors converge to zero as $\Delta t, \Delta \pmb{X} \rightarrow 0$. To ensure that we sufficiently control them for the discretization resolution, we estimate these errors w.r.t. a reference solution computed on a very fine grid. Section~\ref{sec:num_experiments} reports the values of the errors in the proposed approach.
\section{Numerical Experiments and Results}\label{sec:num_experiments}
This section focuses on the example of the Uruguayan power grid and presents the numerical results using the proposed approach to solve the OC problem for the related coupled power system.
\subsection{Description of the Uruguayan Coupled Power System}
\label{sec:Description of the Uruguyan Coupled Power System}
The considered system is illustrated by Figure~\ref{DamsConnection}. It consists of (i) four dams: a cascade of three connected dams located on the same river (Bonete, Baygorria, and Palmar) and an independent dam (Salto Grande), (ii) four FFSs, and (iii) a single battery.  
\begin{figure}[h!]
	\centering
	\includegraphics[width=1\textwidth]{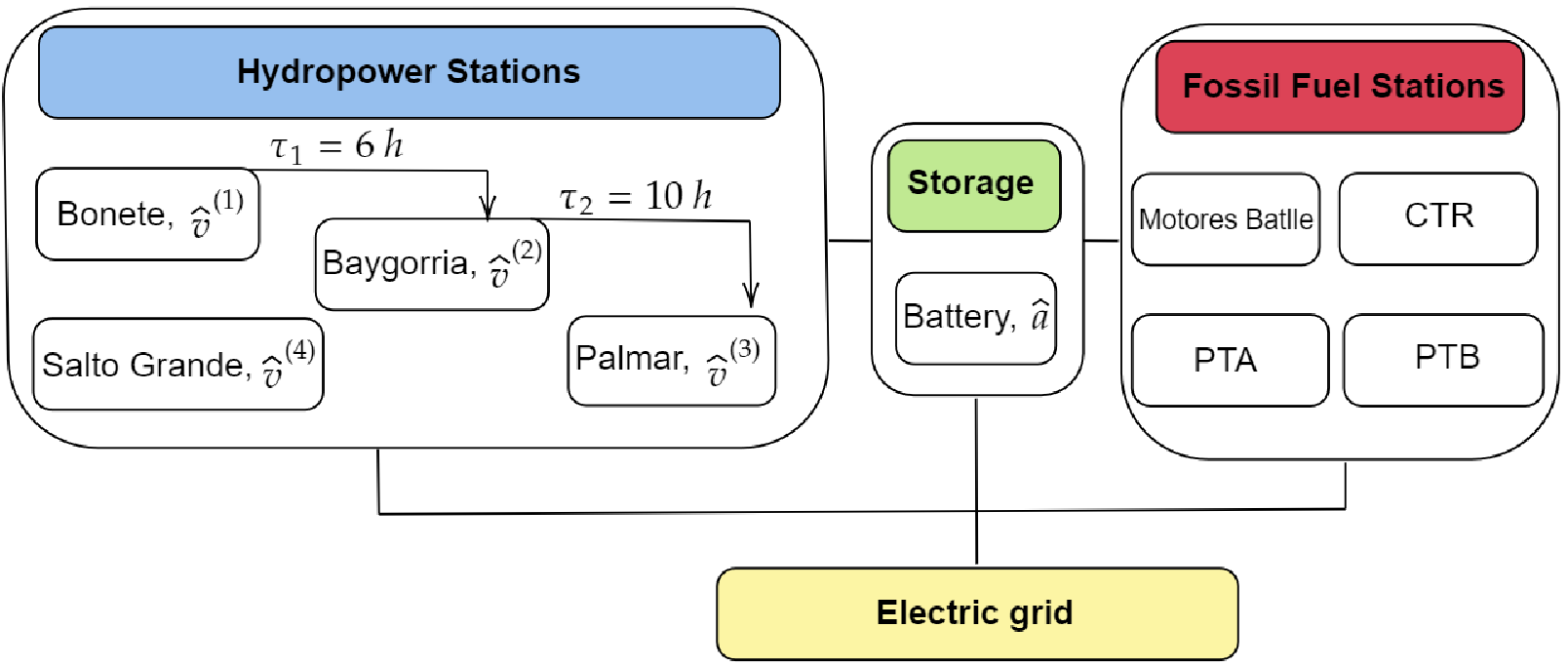}
	\caption{\label{DamsConnection} Network of dams and fossil fuel power stations (FFSs) in the Uruguayan power system.}
\end{figure}
We model the dam dynamics, as discussed in Section~\ref{sec:Dams Modeling}. The time delay between Bonete and Baygorria is $\tau_1 = 6$~h, and the delay between Baygorria and Palmar is $\tau_2=10$~h. The maximum turbine flow and spillage as functions of volume were provided for input by the ADME company (\url{https://adme.com.uy}) and are defined as follows (see Figure~\ref{MaxFlow} for illustration):
\begin{small}
	\begin{align}\label{eq:model_phimax}
		&\bar{\phi}^{(1)}_{Tur}(\hat{v}^{(1)}) = c_2^{(1)}(H^{(1)}(\hat{v}^{(1)})-h_{01})^2 + c_1^{(1)}(H^{(1)}(\hat{v}^{(1)})-h_{01}) + c_0^{(1)}, \nonumber\\
		&\bar{\phi}^{(2)}_{Tur}(\hat{v}^{(2)}) =\begin{cases} c_1^{(2)}(H^{(2)}(\hat{v}^{(2)})-h_{02}) + c_0^{(2)},~~~\text{if }H^{(2)}(\hat{v}^{(2)})-h_{02} \leq 14 \\
			g_1^{(2)}(H^{(2)}(\hat{v}^{(2)})-h_{02}) + g_0^{(2)},~~~\text{if }H^{(2)}(\hat{v}^{(2)})-h_{02} \geq 14
		\end{cases} \nonumber\\ 
		&\bar{\phi}^{(3)}_{Tur}(\hat{v}^{(3)}) =\begin{cases} c_1^{(3)}(H^{(3)}(\hat{v}^{(3)})-h_{03}) + c_0^{(3)},~~~\text{if }H^{(3)}(\hat{v}^{(3)})-h_{03} \leq 23 \\
			g_2^{(3)}(H^{(3)}(\hat{v}^{(3)})-h_{03})^2+g_1^{(3)}(H^{(3)}(\hat{v}^{(3)})-h_{03}) + g_0^{(3)},~~~\text{if }H^{(3)}(\hat{v}^{(3)})-h_{03} \geq 23
		\end{cases} \nonumber\\  
		&\bar{\phi}^{(4)}_{Tur}(\hat{v}^{(4)}) = 4410, \nonumber\\ 
		&\bar{\phi}^{(i)}_S(\hat{v}^{(i)}) = \bar{\phi}_{\text{max}}^{(i)} - \bar{\phi}^{(i)}_{Tur}(\hat{v}^{(i)}),~~~1 \leq i \leq 4,
	\end{align}
\end{small}
where $\overline{\phi}^{(i)}_{\text{max}} $ is the total maximum flow of the $i$th dam, which is given.}

Motivated by the empirical observations, we approximate $H^{(i)}(\hat{v}^{(i)})$ in the model~\eqref{eq:HP_linearized model}, locally, using a second-order polynomial, 
\begin{equation}\label{eq:hight_volume_dependence}
H^{(i)}(\hat{v}^{(i)}) \approx b_2^{(i)}(v^{(i)})^2+b_1^{(i)} v^{(i)}+ b_0^{(i)},
\end{equation}
where $b_2^{(i)}, b_1^{(i)}$, and $b_0^{(i)}$ are determined from the data.
Table~\ref{NE_ModelParams_Dams} lists the
constants $\{\{b_j^{(i)},c_j^{(i)},g_j^{(i)}\}_{j=0}^2, h_{0i}\}_{i=1}^3$. 
\begin{figure}[h!]
\centering
\includegraphics[width=1\textwidth]{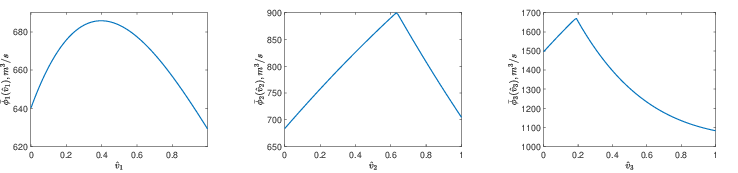}
\caption{\label{MaxFlow}  Maximum turbine flow as a function of volume for the three dams.}
\end{figure} 
The additional constants $\{\underline{v}^{(i)}, \bar{v}^{(i)}\}_{i=1}^4 $ for the minimum and maximum values of the dam volumes, $\{\eta^{(i)},d^{(i)}\}_{i=1}^4$ applied in the power models~\eqref{eq:HP_linearized model} and the water costs $\{K_H^{(i)}\}_{i=1}^4$ are also provided in Table~\ref{NE_ModelParams_Dams}.

The FFS and battery dynamics are modeled as in Sections~\ref{sec:Fossil Fuel Stations (FFS) Modeling} and~\ref{sec:Battery Modeling}, respectively. Moreover, the parameters $\{\bar{P}_F^{(i)}, K_F^{(i)}\}_{i=1}^4$ of the four FFSs (Motores Battle, PTA, PTB, and CTR), and the parameters $\bar{A}, \underline{P}_A$ of the battery are provided in Tables~\ref{NE_ModelParams_FFS} and~\ref{NE_ModelParams_Battery}, respectively.

\begin{figure}[h!]
\centering
\includegraphics[width=0.6\textwidth]{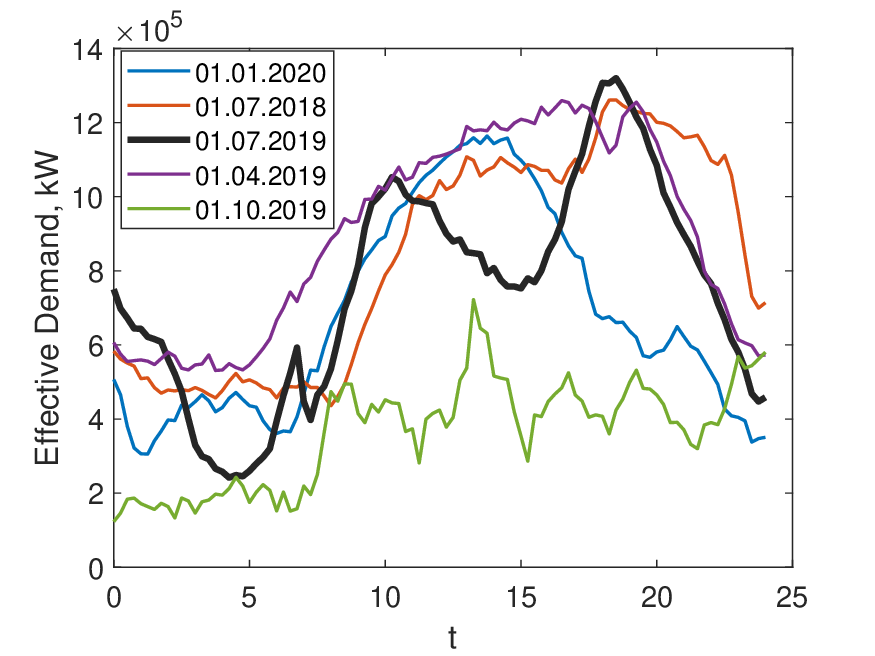}
\caption{\label{NE_Demand} Demand data for different days from all seasons. The results are displayed for January~7, 2019.}
\end{figure}

The effective demand curve in~\eqref{eq:demand_eq} is obtained using linear interpolation from input data (see Figure~\ref{NE_Demand}) containing measurements with a resolution of 10~min. 
We present the numerical results for the data corresponding to January~7, 2019. However, we also tested the proposed approach for various demand data for several different days from all seasons of the year with similar results.

\subsection{Results of the Optimal Control Problem Related to the Uruguayan Coupled Power System}
This section presents the results of solving the OC problem corresponding to the management of the Uruguayan coupled power system described in Section~\ref{sec:Description of the Uruguyan Coupled Power System}. 
We display the results for the various stages of the proposed approach, as described in Sections~\ref{sec:SOC_formulation} and ~\ref{sec:SOC_num_approach}. Section~\ref{sec:Solution of the HJB} explores the behavior of the value function obtained after solving the HJB equation~\eqref{eq:HJB_PDE} for fixed continuous-time Lagrangian multipliers. Section~\ref{sec:Dual optimization results} reveals the optimization results of the dual problem and highlights the effects of the optimal continuous-time Lagrangian multipliers on the system in terms of dual cost and dual gap, also demonstrating the advantage of level refinement of these multipliers in reducing the relative dual gap to the desired gap (i.e., below 2\%). 
Section~\ref{sec:Penalization constant tuning} details the outputs of the penalization procedure performed to smooth the nearly optimal primal controls of interest. Finally, Section~\ref{sec:Optimal solution results} presents the synthetic results of the dual and primal problems in terms of the obtained OCs, power-generation profile over the optimal path, and other related outputs.

\subsubsection{{Value Function of the Dual Problem}}\label{sec:Solution of the HJB}

The numerical solution of HJB~\eqref{eq:HJB_PDE}, 
related to the dual function~\eqref{eq:dual function (markovian problem)}, is obtained with the scheme described in Section~\ref{sec: Monotone Scheme for Solving the HJB}. We use the uniform rectangular grid with $\Delta t = 0.25$~h, $\{\Delta \hat{v}^{(i)}\}_{i=1}^4 = \Delta \hat{a} = 0.25$ (in normalized scale), leading to a Courant number below 0.8 in the CFL condition~\eqref{eq: CFL}. The battery has the fastest dynamics; therefore, the largest contribution to the CFL value is made by the battery. Thus, the choice of $\Delta t$ is almost entirely determined by $\Delta \hat{a}$. Table~\ref{tab: CFL} displays CFL contributions from all variables. This choice of coarse discretization is sufficient to reach the desired error bound~\eqref{eq:error_bound_dual_gap} for a relative duality gap below 2\%. Table~\ref{tab: erros_our_approach} summarizes the values of the errors of the proposed approach.
\begin{table}[h!]
\centering
\begin{tabular}
	{|p{2.8cm}||p{2.8cm}||p{2.8cm}||p{3cm}||p{2.6cm}|}
	\hline
	Bonete, $\hat{v}^{(1)}$  & Baygorria, $\hat{v}^{(2)}$  & Palmar, $\hat{v}^{(3)}$   & Salto Grande, $\hat{v}^{(4)}$ & Battery, $\hat{a}$ \\
	\hline  
	$5.6 \cdot 10^{-4}$ & $3.1 \cdot 10^{-2}$   & $5.5 \cdot 10^{-3}$ &   $2.1 \cdot 10^{-2}$ & $7.1 \cdot 10^{-1}$\\
	\hline  
\end{tabular}
\caption{\label{tab: CFL} Terms $\frac{||f_k||_{\infty}}{\Delta x^{(k)}}$, comprising the Courant--Friedrichs--Lewy (CFL) condition according to~\eqref{eq: CFL} for each $k=1,...,5$ spatial dimension.}
\end{table}

\begin{figure}[h!]
\centering
\includegraphics[width=1\textwidth]{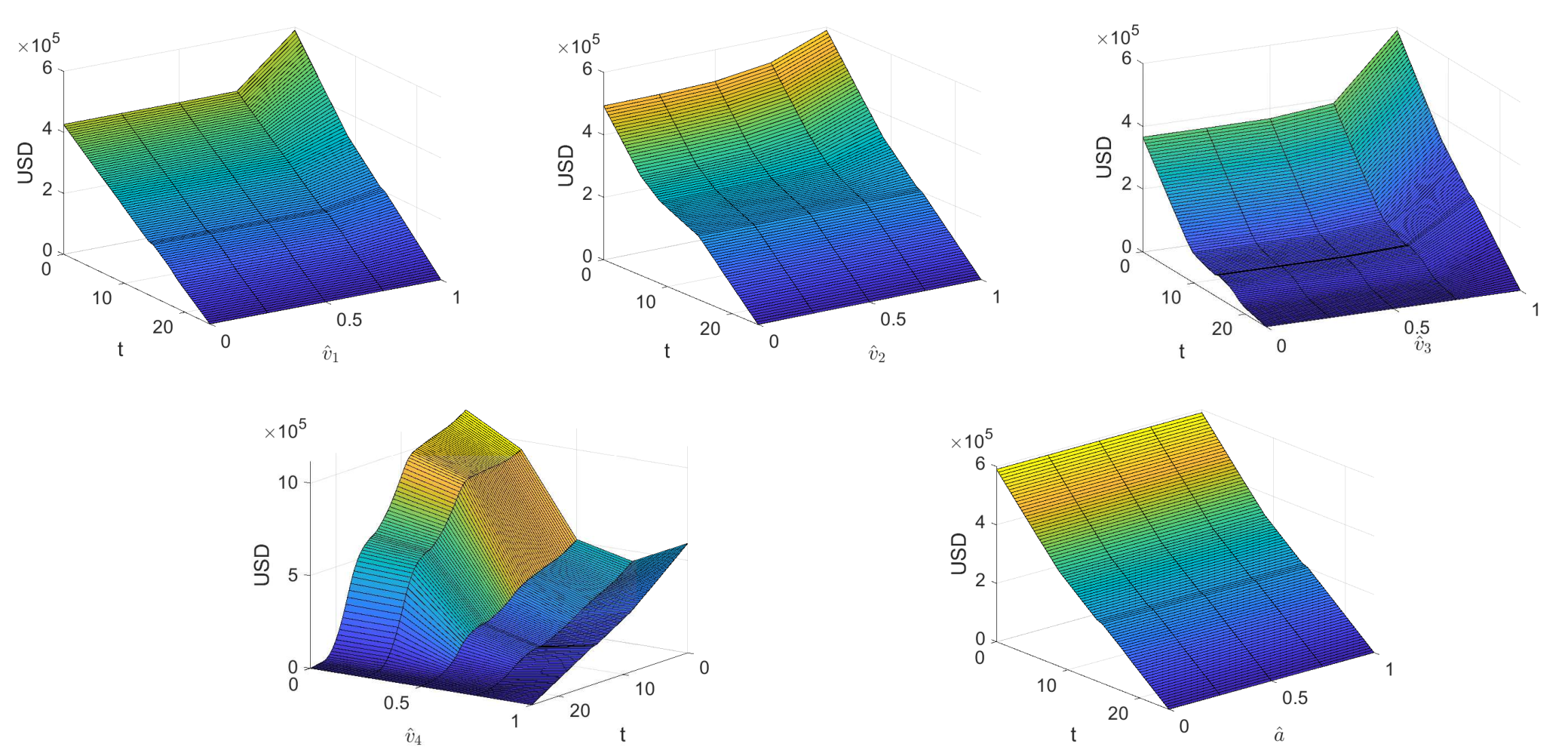}
\caption{\label{NE_ValueFunction} Projections of the value function onto the axes of the time and state variables. For all projections onto the $j$th component of $\bf{X}$, the rest of the variables are fixed at $\{\bf{X}_i\}_{i=1,i\neq j}^5 = 1$.}
\end{figure}

Figure~\ref{NE_ValueFunction} illustrates five slices of the value function with the optimal Lagrange multipliers, $\widetilde{\pmb{\lambda}}^{*}_{(2)}$,  each with four fixed-state variables. 
The observation that the dependence of the value function on $\hat{v}^{(4)}$ seems to be more complex than the dependence on the other state variables 
suggests that we can reduce the dimension of the state space to just one dimension using an ansatz function for the dependence of the value function w.r.t other dimensions and learn its parameters with the available data produced with this coarse discretization grid. This dimension reduction approach is left for future work.

\subsubsection{Dual Optimization Results}
\label{sec:Dual optimization results}
In this context, the Lagrangian multiplier $\boldsymbol{\lambda}(t)$, associated with the relaxation of Constraints~\ref{virtual_controls_target}, 
has two components: $\lambda_{1}(t)$ (related to the water coming from Bonete to Baygorria), and $\lambda_{2}(t)$ (related to the water coming from Baygorria to Palmar). Both $\lambda_{1}(t)$ and $\lambda_{2}(t)$ are approximated by piecewise constant functions, as defined in~\eqref{eq: lambda_approx_first}. 

The goal is to determine the solution of the dual problem that ensures a desired relative dual gap below $\leq 2$\%. The result was achieved using Algorithm~\ref{algo_gen} with two levels of refinement. The optimal configuration in each level is given in Table~\ref{NE_Lam_opt}, and the corresponding values of the dual and primal costs and the relative dual gap are summarized in Table~\ref{NE_DGap_opt}.
\begin{table}[h!]
\centering
\begin{tabular}
	{ |p{4cm}||p{10cm}|}
	\hline
	Parameters & Value\\
	\hline  
	\multicolumn{2}{|c|}{Refinement level $n=1$}\\
	\hline 
	$m_1$, $m_2$ & 1 \\
	$N_{\text{iter}}$ & 30  \\
	$\widetilde{\pmb{\lambda}}^{(0)}_{(1)}$ &  $(10^{-4}, 10^{-4})$\\
	$\widetilde{\pmb{\lambda}}^{*}_{(1)}$ &  $(-1.68 \times 10^{-4}, -8.92 \times 10^{-4})$\\
	\hline  
	\multicolumn{2}{|c|}{Refinement level $n=2$}\\
	\hline 
	$m_1$, $m_2$ & 2 \\
	$N_{\text{iter}}$ & 30 \\
	$\widetilde{\pmb{\lambda}}^{(0)}_{(2)}$ & $((-1.68 \times 10^{-4}, -1.68 \times 10^{-4}),(-8.92 \times 10^{-4}, -8.92 \times 10^{-4}))$  \\
	$\widetilde{\pmb{\lambda}}^{*}_{(2)}$ & $((-2.44 \times 10^{-3}, -3.89 \times 10^{-6}),(-2.16 \times 10^{-3}, -2.35 \times 10^{-3}))$  \\
	\hline  
\end{tabular}
\caption{\label{NE_Lam_opt} Limited memory bundle method (LMBM) optimization outputs for each refinement level.}
\end{table}
\begin{table}[h!]
\centering
\begin{tabular}
	{ |p{3.2cm}||p{3.2cm}||p{3.5cm}||p{3cm}|}
	\hline
	Refinement level &  Dual cost  {(in USD)} & Primal cost {(in USD)} & Relative dual gap\\
	\hline  
	$n=1$ & $5.4559 \times 10^5$  & $5.9555 \times 10^5$ & $9.16$\% \\
	$n=2$ & $5.8952 \times 10^5$  & $5.9555 \times 10^5$ & $1.02$\%  \\
	\hline  
\end{tabular}
\caption{\label{NE_DGap_opt} Dual gap reduction using dual problem optimization with level refinement.}
\end{table}

\begin{figure}[h!]
\centering
\includegraphics[width=0.9\textwidth]{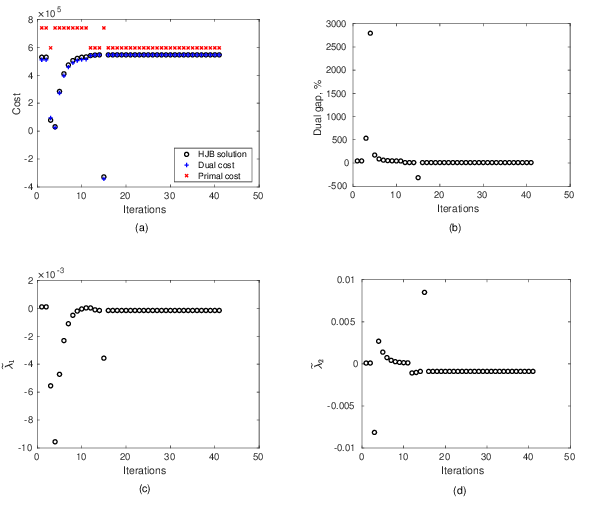}
\caption{\label{NE_LamOpt_R} Dual optimization results for the refinement level at $n=1$: (a) the HJB solution/dual cost/primal cost (smooth), (b) relative dual gap, and (c), (d) components of $\widetilde{\pmb{\lambda}}_{(1)} \in \mathbbm{R}^2$ w.r.t. iterations.}
\end{figure}

Figure~\ref{NE_LamOpt_R} depicts the HJB solution, dual gap, and components of $\widetilde{\pmb{\lambda}}_{(1)}$ at each iteration of the LMBM optimization algorithm. As observed in the plots, not every step is taken in the direction of maximization. At iteration $k=15$, the dual cost $\widetilde{\theta}(\widetilde{\pmb{\lambda}}_{(1)}^{(k)})$ jumps away from the optimum. This behavior illustrates the difficulty of the nonsmooth optimization (Problem~\ref{prob:dual}). The LMBM can quickly recover from this dip with the accumulated information about the subgradient on previous iterations. 

\begin{figure}[h!]
\centering	
	\includegraphics[width=0.9\textwidth]{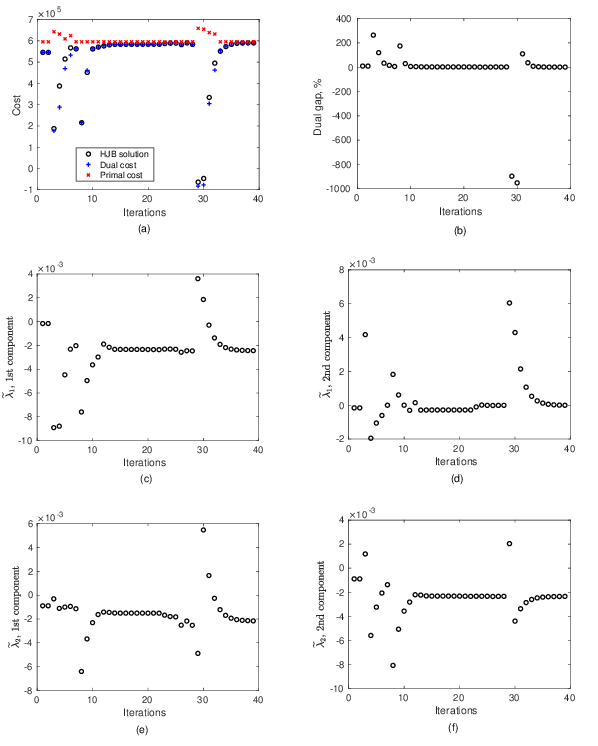}
\caption{\label{NE_LamOpt_R2} Dual optimization results for the refinement level at $n=2$: (a) the HJB solution/dual cost/primal cost (smooth), (b) relative dual gap, and (c), (d), (e), (f) components of $\widetilde{\pmb{\lambda}}_{(2)} \in \mathbbm{R}^4$ w.r.t. iterations.}
\end{figure}

Figure~\ref{NE_LamOpt_R2} illustrates plots similar to those in Figure~\ref{NE_LamOpt_R}, but for the next refinement level where $\widetilde{\theta}(\widetilde{\pmb{\lambda}}_{(2)})$ is optimized in four dimensions. We can verify that the optimization becomes more challenging because taking steps in the wrong direction occurs more frequently than in the previous case. The LMBM method can still effectively address the jumps and obtain the desired dual gap reduction due to the accumulated information about the subgradient.

\subsubsection{Results of the Penalization Parameters Tuning Procedure}\label{sec:Penalization constant tuning}

\begin{figure}[h!]
\centering
	\includegraphics[width=1\textwidth]{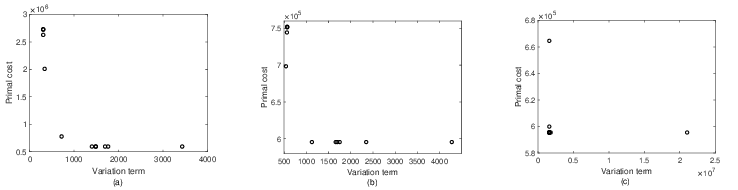}
\caption{\label{NE_LCurves} The L-curve plots corresponding to the optimal $\widetilde{\pmb{\lambda}}^{*}_{(2)}$ for (a) turbine flow, (b) spillage and (c) the battery controls. The optimal penalizing constants are $\beta_{Tur}=10$, $\beta_S=10$, and $\beta_A=10^4$, which correspond to the point in the corner of the L-curve in each plot.}
\end{figure}

This section presents the results of the penalization parameter tuning described in Section~\ref{sec:Admissible Controls Smoothing: The Penalization Parameter Tuning Procedure}. Due to the different nature of the controls (turbine, spillage, and battery), the tuning procedure is performed separately for the three groups of controls with coefficients $\beta_{Tur}^{(i)}=\beta_{Tur},~i=1,...,4$ for the turbine controls, $\beta_S^{(i)}=\beta_S,~i=1,...,4$ for the spillage controls, and $\beta_A$ for the battery control. As in Algorithm~\ref{algo_gen}, this tuning procedure is conducted at each refinement level, and the corresponding outputs are provided in Table~\ref{NE_beta_tuning}. The L-curves produced at the refinement level of $n=2$ are depicted in Figure~\ref{NE_LCurves}, suggesting the optimal penalization constants are $\beta_{Tur}=10$, $\beta_S=10$, and $\beta_A=10^4$. We observe that $\beta_A>> \beta_{Tur}, \beta_S$ due to the high variability of the battery control solutions, as displayed in Figure~\ref{NE_SmoothBat}.
\begin{table}[h!]
\centering
\begin{tabular}
	{ |p{5cm}||p{3cm}||p{3cm}||p{3cm}|}
	\hline
	& Turbine (in USD/($m^6/s^3$) )& Spillage (in USD/($m^6/s^3$) ) & Battery (in USD/($kJ^2/s^3$) )\\
	\hline 
	Initial $\widetilde{\pmb{\lambda}}_0$, $n=0$ & $\beta_{Tur} = 10$ & $\beta_S = 10$ & $\beta_A = 10^2$ \\
	Refinement level $n=1$ & $\beta_{Tur} = 10$ & $\beta_S = 1$ & $\beta_A = 10^4$ \\
	Refinement level $n=2$ & $\beta_{Tur} = 10$ & $\beta_S = 10$ & $\beta_A = 10^4$ \\
	\hline  
\end{tabular}
\caption{\label{NE_beta_tuning} Penalization coefficients $\pmb{\beta}:=(\beta_{Tur}^{(1)},\dots, \beta_{Tur}^{(4)}, \beta_S^{(1)},\dots, \beta_S^{(4)},\beta_A )$ at each refinement level.}
\end{table}
\subsubsection{Optimal Solution}
\label{sec:Optimal solution results}
With the optimal Lagrangian multipliers $\widetilde{\pmb{\lambda}}^{*}_{(2)}$ reported in Table~\ref{NE_Lam_opt}, the dual cost is $5.8952 \times 10^5$~USD, whereas the primal cost is $5.9555 \times 10^5$, implying a relative dual gap of 1.02\% and a total error bound~\eqref{eq:error_bound_dual_gap} for the duality gap below 2\%, with details of each error in~\eqref{eq:error_bound_dual_gap} in Table~\ref{tab: erros_our_approach}.
\begin{table}[h!]
\centering
\begin{tabular}
	{ |p{3.5cm}| |p{3.5cm}||p{2.2cm}||p{2.2cm}||p{2.2cm}|}
	\hline
	Total relative dual gap & Relative dual gap (Error(II)) & Error (I) & Error (III) & Error (IV)\\
	\hline  
	$1.6$\%& $1.02$\%& $5\times 10^{-4}$ \%  & $0.4$\% &   $0.18$\% \\
	\hline  
\end{tabular}
\caption{\label{tab: erros_our_approach} Values of the terms comprising the bound on the dual gap in~\ref{eq:error_bound_dual_gap}.}
\end{table}

\begin{figure}[h!]
\centering
	\includegraphics[width=1\textwidth]{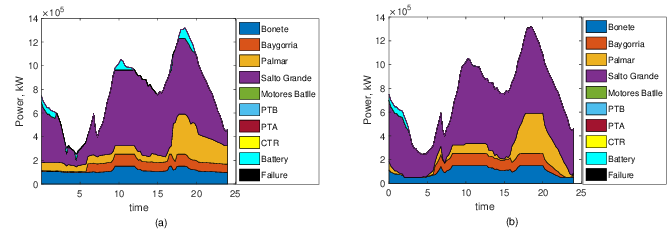}
\caption{\label{NE_Power} Demand and power-generation profile over the optimal (a) and smoothed (b) path. In each case, the demand is satisfied entirely with the power generated by the dams and the battery, for this reason no FFS is activated.}
\label{fig:demand and power curves_nr_2}
\end{figure}

\begin{figure}[h!]
\centering
		\includegraphics[width=1\textwidth]{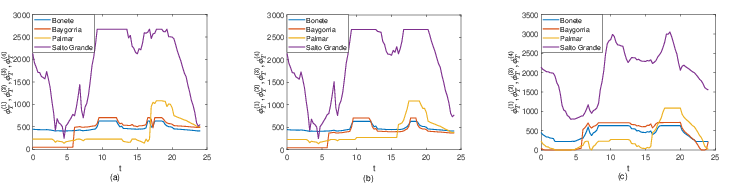}
\caption{\label{NE_SmoothFlow} Turbine flow solutions: (a) dual optimal, (b) primal admissible, and (c) primal smoothed.}
\end{figure}

\begin{figure}[h!]
\centering
	\includegraphics[width=1\textwidth]{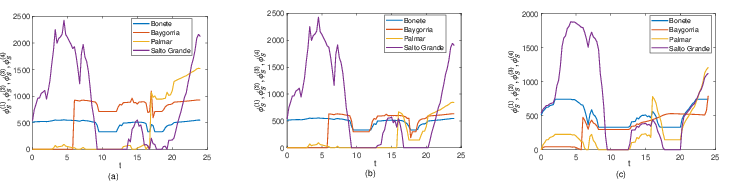}
\caption{\label{NE_SmoothSpil} Spillage solutions: (a) dual optimal, (b) primal admissible, and (c) primal smoothed.}
\end{figure}

\begin{figure}[h!]
\centering
	\includegraphics[width=1\textwidth]{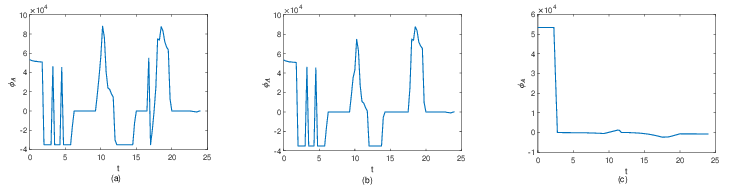}
\caption{\label{NE_SmoothBat} Battery control solutions: (a) dual optimal, (b) primal admissible, and (c) primal smoothed. Smoothing has especially strong impact on the battery controls: plot (c) shows a better practical solution while maintaining a small relative dual gap.}
\end{figure}

\begin{figure}[h!]
\centering
	\includegraphics[width=0.8\textwidth]{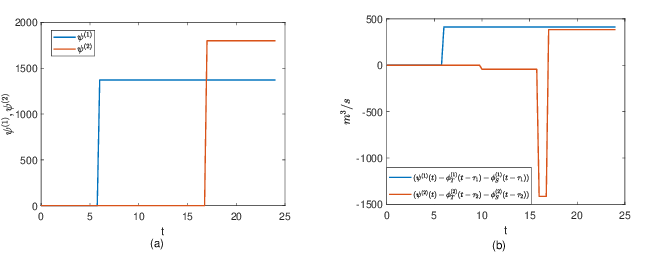}
\caption{\label{NE_VirtCont} (a) Virtual controls. (b) Difference between virtual and real water. }
\end{figure}

\begin{figure}[h!]
\centering
	\includegraphics[width=0.8\textwidth]{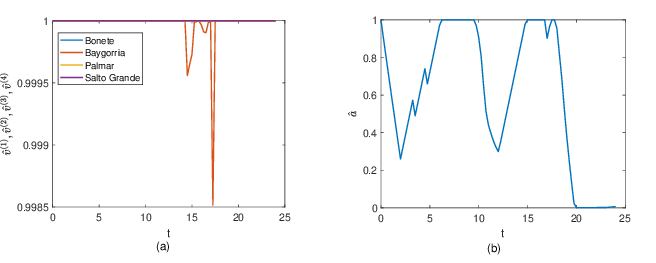}
\caption{\label{NE_States} State variables over the optimal path: (a) dam volumes and (b) battery charge.}
\end{figure}

The outputs of the optimal solution are illustrated in Figures~\eqref{NE_Power} to~\eqref{NE_States}, including (i) the OCs related to the turbine flow, spillage, FFSs, and battery, besides the virtual controls; (ii) the volumes of the dams and battery charge; and (iii) the profile of the power production over the optimal and smoothed path.

From the power profile in Figure~\ref{NE_Power}a, we conclude that most of the power is generated by the largest dam (Salto Grande). The most expensive power generators (FFSs) are unnecessary because the demand can be satisfied with the dams and battery only. The battery is primarily used in three intervals: at the beginning, around 10 h, and around 18 h, corresponding to the peaks of electricity demand.

Figures~\ref{NE_SmoothFlow} and~\ref{NE_SmoothSpil} display the turbine flow and spillage controls. The highest flow is associated with the largest dam: Salto Grande. Baygorria is the cheapest dam; however, it is not used much in the first 14 h because, at the initial time, the dams are full, creating the need to turbine and spill water from the dams to prevent them from overflowing. We compared the rates of turbine flow and spillage with the inflow at each dam: Bonete: 958~m$^3$/s, Baygorria: 43~m$^3$/s, Palmar: 226~m$^3$/s, Salto Grande: 2675~m$^3$/s. All decisions about turbining and spilling water in this work are driven by the need to satisfy the corresponding constraint on the volume of the dams rather than the optimality of the solution. 
The turbine and spillage flows in Figures~\ref{NE_SmoothFlow}(a) and~\ref{NE_SmoothSpil}(a) correspond to the optimal solution of the dual problem and, therefore, may be outside of the admissible set due to the relaxation of Constraints~\ref{virtual_controls_target}. Thus, the admissible turbine and spillage flow controls (see Figures~\ref{NE_SmoothFlow}(b) and ~\ref{NE_SmoothSpil}(b)) are computed as a (nearly optimal) solution of the primal problem using the procedure described in Section~\ref{sec:Construction of the admissible Controls for the Primal Problem}. These admissible controls exhibit many jumps and irregularities; therefore, we smooth these controls as explained in Sections~\ref{sec:Admissible Controls Smoothing: Computation of Penalized Controls} and~\ref{sec:Admissible Controls Smoothing: The Penalization Parameter Tuning Procedure} to obtain more convenient solutions in practice, as presented in Figures~\ref{NE_SmoothFlow}(c) and~\ref{NE_SmoothSpil}(c).

The most significant effect of the smoothing was observed on the battery controls (Figure~\ref{NE_SmoothBat}). 
The dual optimal and admissible controls are highly oscillating, making them impractical because the operator must switch the battery from discharging to charging mode too frequently. The smoothed control is a better practical solution without sacrificing too much cost, with a relative dual gap of only 1.02\%.

Figure~\ref{NE_VirtCont} illustrates how the virtual controls work. Consider the connection between Bonete and Baygorria and the virtual water $\psi^{(1)}$. The Hamiltonian is linear w.r.t. the virtual controls, and whether we aim to activate the virtual control is determined by the sign of the term $(\partial_{\hat{v}^{(2)}}u(t) + \lambda_{1}(t+\tau_1))$ in the Hamiltonian (\ref{eq:HJB_PDE}). If this term is negative, $(\partial_{\hat{v}^{(2)}}u(t) + \lambda_{1}(t+\tau_1)) < 0$, the water from the upstream dam is valuable. It takes 6~h for the water from Bonete to reach Baygorria, and $(\partial_{\hat{v}^{(2)}}u(t) + \lambda_{1}(t+\tau_1)) < 0$, for all $t \in [6,24]$~h; thus, we activate the control in this time interval. Now consider the link between Baygorria and Palmar and the virtual water $\psi^{(2)}$. The term $(\partial_{\hat{v}^{(3)}}u(t) + \lambda_{2}(t+\tau_2)) < 0$ in the interval of $t \in [10,24]$~h; thus, we activated the virtual control $\psi^{(2)}$ on this interval. However, the dam is full, and we can only receive more water when we start turbining more water with Palmar due to the high demand starting at around 16~h. The difference between the virtual and real water is presented in Figure~\ref{NE_VirtCont}(b). As we optimize the dual function $\theta(\pmb{\lambda}(t))$ w.r.t. $\pmb{\lambda}$, this difference decreases.

Figure~\ref{NE_States} presents the optimal trajectory of the state variables. As observed in Figure~\ref{NE_States}(a), the volumes of Bonete $\hat{v}^{(1)}$, Palmar $\hat{v}^{(3)}$, and Salto Grande $\hat{v}^{(4)}$ stay at maximum levels because the outflow is equal to the inflow. The battery charge plot in Figure~\ref{NE_States}(b) illustrates how the battery is used throughout the optimal path and depletes at the end.

\section{Conclusions}\label{sec:conclusions}
This work modeled a large-scale power system, including a cascade of hydropower stations, FSSs, and a storage unit represented by a single battery. We provided a detailed description of each power-source cost, power generation, and state dynamics. We formulated the primal OC problem for the large-scale power system with time delays and proposed a continuous-time Lagrangian relaxation technique to address the time delays in the hydropower dynamics, resulting in a dual problem formulation. Further, we proposed a heuristic procedure for obtaining sufficiently smooth, nearly optimal admissible controls for the primal problem. We presented the algorithm to the given problem numerically. This algorithm includes the numerical approximation of the HJB equation using the backward-in-time explicit Euler upwind finite-difference method. For the optimization of the dual problem, we applied the LMBM. We showed how to use the L-curve approach to fine-tune the parameters for constructing a sufficiently smoothed primal solution from the dual solution. Finally, we demonstrated the effectiveness of the proposed method by presenting numerical results obtained through its application to the management of a partial model of the Uruguayan power grid. This practical application is an illustrative example of the real-world viability and efficacy of the proposed approach.

A possible next research direction consists in including the stochastic dynamics of the wind and solar power production, as well as the demand, based on It\^o SDEs, which results in a stochastic optimal control (SOC) problem. The relaxation of the primal SOC problem leads to a stochastic optimization problem, which requires utilization of corresponding optimization methods, e.g., stochastic subgradient method. 
Our numerical results suggest a potential to reduce the dimensionality of the HJB equation by proposing a suitable ansatz for the dependence of the value function on all dimensions except one, transforming the PDE into an ODE. Another dimension reduction technique which is useful in our context is the Markovian projection: an approach that allows to project the wind, solar and demand dynamics onto a single stochastic process. These avenues remain open for future investigation.

~

\textbf{Acknowledgments}
This work was supported by the KAUST Office of
Sponsored Research (OSR) under Award No. URF/1/2584-01-01 and the
Alexander von Humboldt Foundation. E.~Rezvanova, E.~von Schwerin, and R.~Tempone are
members of the KAUST SRI Center for Uncertainty Quantification in
Computational Science and Engineering.

\bibliographystyle{plain}
\bibliography{10_SOC_energy_manuscript} 

\appendix

\section{Remarks on Stochastic Modeling}\label{sec:Remarks on Stochastic Modeling}
\begin{remark}[On the stochastic modeling of the demand and exports]\label{rem: Stochastic modeling of demand and exports}
	In a more detailed and complex model, the demand and export predictions are affected by stochastic perturbations, for instance, i) weather effects on electricity consumption and the dynamics of wind and solar power production, and ii) possible bids on the excess of the energy supply allocated to exports. Modeling the stochastic dynamics of both quantities is left for future work.
\end{remark}

\begin{remark}[On stochastic dynamics]\label{rem:v_SDE}
	One method of representing the uncertainty in the natural input of water and the 
	total volume measurement is to add noise to the dynamics, resulting in It\^{o} stochastic differential equations. Then, the ordinary differential equation (ODE)~\eqref{eq: dams_independent_dynamcis} is replaced by 
	\begin{align}
		d\hat{v}^{(i)}(t) &= \frac{1}{\bar{v}^{(i)} -  \underline{v}^{(i)}} \left( \left(I^{(i)}(t)- \phi_{Tur}^{(i)}(t)- \phi_S^{(i)}(t)\right) dt+\sigma_v^{(i)}\left(t,\hat{v}^{(i)}(t)\right)dW_H^{(i)}\right), && 0 \le t \le  T,
	\end{align}
	where $W_H^{(i)}$ is the $i$th component of a standard $N_D$-dimensional Brownian motion, and the diffusion coefficient is $\sigma_v^{(i)}/\left(\bar{v}^{(i)}-\underline{v}^{(i)}\right)$. The requirement that $0\le \hat{v}^{(i)}(t) \le 1$, for $t>0$, is enforced, almost surely, by modeling $\sigma_v^{(i)}(t,0) = \sigma_v^{(i)}(t,1)=0$, $\forall\: t \in [0,T]$, and maintaining the outflow constraints~\eqref{eq: dams_independent_constraints}.
	We emphasize that $\sigma_v^{(i)}$ should belong to a parametric family of functions that can be calibrated against data, such that (s.t.) it vanishes at the boundaries (i.e., $\sigma_v^{(i)}(t,0) = \sigma_v^{(i)}(t,1)=0$, for $0 \le t \le T$).
	
	Analogously, the diffusion term $\sigma_v^{(i)}(t,\hat{v}^{(i)}(t))/\left(\bar{v}^{(i)} - \underline{v}^{(i)}\right)dW_H^{(i)}$ is added to~\eqref{eq: dams_connected_dynamcis} to obtain the stochastic dynamics for the downstream dams. The time delays in the controls lead to a non-Markovian stochastic OC (SOC) problem. The proposed approach to these problems is analogous to the one used in this paper to deal with the time delays in the deterministic dam dynamics. We leave a more detailed description of the stochastic case with corresponding numerical tests for future work.
	
	The stochastic nature of the dam dynamics does not violate the mass conservation law. 
	The approach considers the uncertainty related to the natural input of water (inflow) or measurements
	~\cite{marton2015analysis, kwon2016probabilistic,carpentier2018dam}. There are alternative ways to represent this uncertainty, for instance, by adding noise to the inflow $I$ in~\eqref{eq: dams_independent_dynamcis} and~\eqref{eq: dams_connected_dynamcis}, which may result in a more complex structure of the SOC. The coefficients $\sigma_v^{(i)}$ can also be considered artificial diffusion coefficients that may add more regularity   to the viscosity solution of the associated HJB equation.
\end{remark}

\begin{remark}[A more general model for battery dynamics] Battery dynamics can account for a possible loss of power and the eventual uncertainty in the charge level. In this case, the dynamics can be modeled as
	\begin{equation}\label{eq:general battery_dynamics}
		d\hat{a}(t)=   - \frac{\overline{P}_A}{\bar{A}} \hat{\phi}_A(t) \left(   1- \delta_{A}   \mathbf{1}_{\{ \hat{\phi}_A(t) <0\}}\right) dt + \hat{a}(t) \sigma_{A} \,dW_A, \quad 0 \le t \le  T,
	\end{equation}
	where $\delta_{A}>0$ is a power loss factor, $\sigma_{A}$ denotes the diffusion coefficient of the battery (associated with the capacity uncertainty) s.t. $\sigma_{A}(0)=\sigma_{A}(1)=0$, and $W_A$ is the corresponding Brownian motion.
\end{remark}
\section{Formulation of the Limited Memory Bundle Method Optimization}\label{sec:Formulation of the LMBM Optimization Method}
This method is designed to solve the optimization problem
\begin{gather}
	\underset{\bf{x}}{\min} ~\theta(\bf{x})
\end{gather}
for a nonsmooth and possibly non-convex function $\theta: \mathbb{R}^n \rightarrow \mathbb{R}$.

The method is a hybrid of the variable metric bundle methods and limited memory variable metric methods \cite{LMBM_PhD,haarala2007globally,karmitsa2012comparing,bagirov2014introduction}. Variable metric bundle methods were developed for nonsmooth optimization and are based on minimizing the linearization of the function at each point. A rather expensive quadratic direction-finding problem must be computed at every iteration to minimize such linearization, which is a time-consuming procedure. Furthermore, large matrices containing information about the metric at the current point must be stored. Limited memory variable metric methods avoid these difficulties by computing the search direction using a limited memory approach, but they were constructed for smooth functions. The limited memory bundle method (LMBM) takes the modified approach from the limited memory variable metric methods and applies it to nonsmooth function optimization.

\textbf{Search direction}. The search direction vector in the LMBM method is calculated as follows:
\begin{gather}
	\pmb{d}_k = -D_k \Tilde{\pmb{\xi}}_k,
\end{gather}
where $\Tilde{\pmb{\xi}}_k$ is the aggregate subgradient at $\pmb{x}_k$ ($\Tilde{\pmb{\xi}}_k = \pmb{\xi}_k$ for a serious step), and $D_k$ is the limited memory variable metric update, which for a smooth function, corresponds to the approximation of the inverse of the Hessian matrix. The matrix $D_k$ is not computed explicitly, (for details see \cite{LMBM_PhD}).

\textbf{Step size.} First, the special line search procedure \cite{LMBM_PhD} generates two points to take a step in the direction of $\pmb{d}_k$:
\begin{gather}
	\pmb{x}_{k+1} = \pmb{x}_k + t^k_L \pmb{d}_k,~~~\text{and}~~~
	\pmb{y}_{k+1} = \pmb{x}_k + t^k_R \pmb{d}_k
\end{gather}
with $\pmb{y}_1 = \pmb{x}_1$, where $t^k_R \in (0,t^k_I]$, $t^k_L \in [0, t^k_R]$ are step sizes, and $t^K_I \in [t_{min}, t_{max})$ is the initial step size. Two types of steps can be taken: the serious and null steps. The step is called serious if we move to a different point by setting $\pmb{x}_{k+1} = \pmb{y}_{k+1}$ on the next iteration. The step is called null if, instead, we remain at the same point $\pmb{x}_{k+1} = \pmb{x}_{k}$ but incorporate the new information about the subdifferential into the aggregate subgradient. The serious step is taken if the following conditions are satisfied
\begin{gather}
	t^k_L=t^k_R>0,~~~~\text{and}~~~\theta(\pmb{y}_{k+1}) \leq \theta(\pmb{x}_{k}) - \epsilon^k_Lt^k_Rw_k,
\end{gather}
where $\epsilon^k_L \in (1,1/2)$, and $w_k$ \eqref{eq: lmbm_stop} is the desirable amount of descent. Otherwise, the null step is taken. In addition, $t^k_L$ and $t^k_R$ are computed via a modified line search procedure (for details, see \cite{LMBM_PhD}).

\textbf{Aggregate subgradient.} 
The subgradient aggregation procedure relies on determining the multipliers $\lambda_i^k$, which minimize the function:
\begin{gather}
	\phi(\lambda_1,\lambda_2,\lambda_3) = [\lambda_1\pmb{\xi}_m + \lambda_2\pmb{\xi}_{k+1} + \lambda_3\Tilde{\pmb{\xi}}_k]^T D_k[\lambda_1\pmb{\xi}_m + \lambda_2\pmb{\xi}_{k+1} + \lambda_3\Tilde{\pmb{\xi}}_k] + 2(\lambda_2\beta_{k+1} + \lambda_3\Tilde{\beta}_k),
\end{gather}
where the subgradient locality measure is $\beta_k = \max\{|\theta(\pmb{x}_k) - \theta(\pmb{y}_{k+1}) + \pmb{\xi}^T_{k+1}(\pmb{y}_{k+1} -\pmb{x}_k)|,\gamma||\pmb{y}_{k+1}-\pmb{x}_k||^w\}$. For convex functions, $\gamma = 0$, and $w \geq 1$ is the locality measure parameter determined by the user.
The aggregate subgradient $\Tilde{\pmb{\xi}}_{k+1}$ is a convex combination of the subgradients computed up to the current iteration:
\begin{gather}
	\Tilde{\pmb{\xi}}_{k+1} = \lambda_1^k \pmb{\xi}_m + \lambda_2^k \pmb{\xi}_{k+1} + \lambda_3^k \Tilde{\pmb{\xi}}_{k},
\end{gather}
where $\pmb{\xi}_m$ denotes the current subgradient at $\pmb{x}_k$, $\pmb{\xi}_{k+1}$, and the auxiliary subgradient at $\pmb{y}_{k+1}$ and $\Tilde{\pmb{\xi}}_{k}$ represents the current aggregate subgradient, with $\Tilde{\pmb{\xi}}_1 = \pmb{\xi}_1$.  

\textbf{Stopping criteria.} 
A number of stopping criteria exist for this algorithm. The primary stopping criterion is the following:
\begin{align} \label{eq: lmbm_stop}
	w_k = -\Tilde{\pmb{\xi}}_k^T \pmb{d}_k + 2\Tilde{\beta}_k \leq \epsilon,
\end{align}  
where $\epsilon$ is the user-specified tolerance.
\section{Parameters of the Numerical Problem}\label{sec:Parameters of the Numerical Problem}
\begin{table}[h!]
	\centering
	\begin{tabular}
		{ |p{2.5cm}||p{2.5cm}||p{2.5cm}||p{3.5cm}||p{3.5cm}|}
		\hline
		Parameter &Bonete, $i=1$& Baygorria, $i=2$&Palmar, $i=3$& Salto Grande, $i=4$\\
		\hline  
		$b^{(i)}_2$, $b^{(i)}_1$, $b^{(i)}_0$ & -3.77, 15.7, 69.6  & -1.4, 8.89, 47.5 &   -7.5, 23.7, 25.8 & -21.2, 54.8, 1.92\\
		$c^{(i)}_2$, $c^{(i)}_1$, $c^{(i)}_0$ & -1.4, 61.1, 35.8  & -, 42.9, 300.5 & -, 42.0, 707.7 & -, -, -\\
		$g^{(i)}_2$, $g^{(i)}_1$, $g^{(i)}_0$ & -, -, -  & -, -81.3, 2039.2 & 2.3, -184.6, 4684.4 & -, -, -\\
		$h_{0i}$, m & 53.96  & 38.57 & 7.00 & 5.60\\
		$d^{(i)}$ & 6.8 $\cdot 10^{-4}$ & 14.0 $\cdot 10^{-4}$& 13.0 $\cdot 10^{-4}$& 14.7 $\cdot 10^{-4}$\\
		$\eta^{(i)}$ & 8.86  & 9.35 & 9.30 & 10.23\\
		$\overline{\overline{\phi}}_{i}$, m$^3$/s & 1371.74  & 1799.68  &3344.70 & 8820 \\
		$\overline{v}^{(i)}$, m$^3$ & 10.7 $\cdot 10^9$  & 0.678 $\cdot 10^9$ & 3.53 $\cdot 10^9$ & 5.18 $\cdot 10^9$\\
		$\underline{v}^{(i)}$, m$^3$ & 1.85 $\cdot 10^9$  & 0.47 $\cdot 10^9$ & 1.36 $\cdot 10^9$ & 3.65 $\cdot 10^9$\\
		$K_{i}$, USD/m$^3$ & 12 $\cdot 10^{-4}$  & 3.96 $\cdot 10^{-7}$  & 23 $\cdot 10^{-4}$  & 16 $\cdot 10^{-4}$ \\
		\hline  
	\end{tabular}
	\caption{\label{NE_ModelParams_Dams} Parameters and coefficients in modeling the dams.}
\end{table}

\begin{table}[h!]
	\centering
	\begin{tabular}
		{ |p{4cm}||p{2.5cm}||p{2.5cm}||p{2.5cm}||p{2.5cm}|}
		\hline
		Parameter & Motores Batlle, $i=1$& PTA, $i=2$& PTB, $i=3$& CTR, $i=4$\\
		\hline  
		$\bar{P}_F^{(i)}$, kW & $7 \cdot 10^4$ & $28.8 \cdot 10^4$ & $36 \cdot 10^4$ & $20 \cdot 10^4$ \\
		$K_F^{(i)}$, USD/MWh & 131  & 193.7  & 189.2  & 222.7 \\
		\hline  
	\end{tabular}
	\caption{\label{NE_ModelParams_FFS} Parameters and coefficients in modeling the fossil fuel power station (FFS).}
\end{table}

\begin{table}[h!]
	\centering
	\begin{tabular}
		{ |p{4cm}||p{5cm}||p{5cm}|}
		\hline
		Parameter & $\overline{A}$, kWh & $\underline{P}_A$, kW\\
		\hline 
		Value &  $1.4 \cdot 10^5$ & $10^5$ \\
		\hline  
	\end{tabular}
	\caption{\label{NE_ModelParams_Battery} Parameters and coefficients in modeling the battery.}
\end{table}
\section{Monotone Scheme for HJB}\label{sec:Monotone Scheme for HJB}
	We consider $(t,\mathbf{x}) \in (0,T) \times [0,1]^{N_D+1} $, $\mathbf{x} := (x^{(1)},...,x^{(N_D+1)})$, and introduce the following set of indices in $\nset^{N_D+1}$, $\mathbf{i}:=(i_1,\dots,i_{N_D+1})$, $\mathbf{i}_k^+=(i_1,\dots,i_{k-1},i_k+1,i_{k+1},\dots,i_{N_D+1})$, and $\mathbf{i}_k^-:=(i_1,\dots,i_{k-1},i_k-1,i_{k+1},\dots,i_{N_D+1})$. We let $\{M_k\}_{k=1}^{N_D+1} \in \nset^{N_D+1}_{+}$ and $N \in \nset_{+}$ s.t. $\{\Delta x^{(k)} = 1/M_k\}_{k=1}^{N_D+1}$ and $\Delta t = T/N$. We consider a uniform grid, $\mathcal{D}:=D^{\Delta t} \times D^{\Delta x^{(1)}} \times \dots \times D^{\Delta x^{(N_D+1)}} $, defined by
\begin{align}
	D^{\Delta x^{(k)}}:=\{	0&=x^{(k)}_{0} < x^{(k)}_{1} = x^{(k)}_{0} + \Delta x^{(k)} < \dots< x^{(k)}_{M_k} =1\}, ~~~1 \leq k \leq N_D+1,\\
	D^{\Delta t}:=	\{0&=t_{0} < t_{1} = t_{0} + \Delta t <  \dots < t_{N} =T\}.
\end{align}
We represent the drift terms of Dynamics~\ref{relaxed_dynamics} 
with $\mathbf{f} = (f_1,\dots,f_{N_D+1})$, where, at each $(t_j, \mathbf{x}_{\mathbf{i}})\in \mathcal{D}$, we have the following:
\begin{align}
	f_k(t_j,\mathbf{x}_{\mathbf{i}},\hat{\boldsymbol{\phi}}_{E}) & = 
	\begin{cases} 
		\frac{I^{(k)}(t_j)- \bar{\phi}_T^{(k)}(x^{(k)}_{i_k}) \hat{\phi}_T^{(k)} (t_j)  -   \bar{\phi}_S^{(k)}(x^{(k)}_{i_k}) \hat{\phi}_S^{(k)} (t_j) + \sum_{l\in\mathcal{B}(k)} \bar{\psi}^{(l)} \hat{\psi}^{(l)}(t_j) }{\bar{v}^{(k)} -  \underline{v}^{(k)}}, 
		& k=1,\dots,N_D, \\
		-  \frac{\hat{P}_A(x^{(k)}_{i_k}) }{\bar{A}} \hat{\phi}_A(t_j), & k = N_D + 1. 
	\end{cases}
\end{align}
Moreover, we define $f^{j^+}_{k,\mathbf{i}}(\hat{\boldsymbol{\phi}}_{E}) := \max\{f_k(t_j,\mathbf{x}_{\mathbf{i}},\hat{\boldsymbol{\phi}}_{E}),0\}$ and $ f^{j^-}_{k,\mathbf{i}}(\hat{\boldsymbol{\phi}}_{E}) := \min\{f_k(t_j,\mathbf{x}_{\mathbf{i}},\hat{\boldsymbol{\phi}}_{E}),0\}$.

We let $U^j_{\mathbf{i}}$ denote the numerical approximation of $u(t,\mathbf{x};\boldsymbol{\lambda})$ in \eqref{eq:HJB_PDE} at point $(t_j,\mathbf{x}_{\mathbf{i}})$. Then, the upwind finite-difference scheme is given by 
\begin{align}\label{eq:upwind_FD_scheme}
	\begin{cases}
		\frac{U^{j-1}_{\mathbf{i}}- U^j_{\mathbf{i}}}{\Delta t} &= \sum_{k=1}^{N_D+1} \Big(f^{{j}^+}_{k,\mathbf{i}}(\hat{\boldsymbol{\phi}}_{E,\mathbf{i}}^j)D^+_{k,\mathbf{i}}U^j + f^{{j}^-}_{k,\mathbf{i}}(\hat{\boldsymbol{\phi}}_{E,\mathbf{i}}^j)D^-_{k,\mathbf{i}}U^j\Big) + L(t_j,\mathbf{x}_{\mathbf{i}},\hat{\boldsymbol{\phi}}_{E,\mathbf{i}}^j;\boldsymbol{\lambda}),\\
		\hat{\boldsymbol{\phi}}_{E,\mathbf{i}}^j &:= \underset{\hat{\boldsymbol{\phi}}_{E} \in \mathcal{A}_{E}(t, \mathbf{X}(t))}{\text{arg}\min} \Big[ \sum_{k=1}^{N_D+1} \Big(f^{{j}^+}_{k,\mathbf{i}}(\hat{\boldsymbol{\phi}}_{E})D^+_{k,\mathbf{i}}U^j + f^{{j}^-}_{k,\mathbf{i}}(\hat{\boldsymbol{\phi}}_{E})D^-_{k,\mathbf{i}}U^j\Big) + L(t_j,\mathbf{x}_{\mathbf{i}},\hat{\boldsymbol{\phi}}_{E};\boldsymbol{\lambda})\Big],
	\end{cases}
\end{align}
where $D^+_{k,\mathbf{i}}U^j $ and $D^-_{k,\mathbf{i}}U^j$ are the upwind derivatives given by
\begin{equation}
	D^+_{k,\mathbf{i}}U^j = \frac{U^j_{\mathbf{i}_k^+} - U^j_{\mathbf{i}}}{\Delta x^{(k)}},~~~~~D^-_{k,\mathbf{i}}U^j = \frac{U^j_{\mathbf{i}} - U^j_{\mathbf{i}_k^-}}{\Delta x^{(k)}}.
\end{equation}
Additionally, $L(t_j,\mathbf{x}_{\mathbf{i}},\hat{\boldsymbol{\phi}}_{E};\boldsymbol{\lambda})$ is the running cost part of~\eqref{eq:objective_fun} evaluated at $(t_j,\mathbf{x}_{\mathbf{i}})$.

To ensure the stability of the upwind finite-difference scheme \eqref{eq:upwind_FD_scheme}, we impose the Courant--Friedrichs--Lewy (CFL) condition 
\begin{equation}\label{eq: CFL}
	\Delta t \sum_{k=1}^{N_D+1} \frac{||f_k||_{\infty}}{\Delta x^{(k)}} \leq 1.
\end{equation}
\begin{remark}
	In general, an upwind finite-difference scheme requires the computational domain to be hyperpyramidal to preserve the upwind structure of the derivatives at the boundaries of the domain. However, in this case, all drift terms $f_k$ naturally have a sign at the boundaries
	allowing the use of the interior point for the derivative approximation.
\end{remark}

\end{document}